\documentclass[11pt,a4paper,headinclude,footinclude,fleqn,reqno]{amsart}                 
\usepackage[T1]{fontenc}                   
\usepackage[utf8]{inputenc}                 
\usepackage[english]{babel}       
\usepackage{graphicx}                      %
\usepackage[font=small]{quoting}            %
\usepackage{caption}  
     \usepackage{amsmath}

 \usepackage[top=.8in,bottom=1.2in,left=1.5in,right=1.5in]{geometry}
\usepackage{verbatim}
\usepackage{color}
\usepackage{picture}
\usepackage{graphicx}
\usepackage{graphics}
\usepackage{comment}
\usepackage{hyperref}
\hypersetup{colorlinks,linkcolor={blue},citecolor={blue},urlcolor={red}}

\newtheorem{thm}{Theorem}[section] 
 
\newtheorem{cor}[thm]{Corollary} 
\newtheorem{prop}[thm]{Proposition}

\theoremstyle{definition} 
\newtheorem{rem}[thm]{Remark} 
\theoremstyle{remark}

\def\e{\epsilon}
\def\S{\Sigma} 
\def\n{\nabla}

\def\p{\partial}

\def\bn{\overline \nabla}

\def\a{\alpha}
\def\b{\beta}
\def\n{\nabla}

\def\o{\omega}
\def\O{\widehat\Sigma}
\def\p{\partial}
\def\e{\epsilon}
\def\a{\alpha}
\def\b{\beta}
\def\g{\gamma}

\def\k{\kappa}

\def\s{\sigma}

\def\n{\nabla}
\def\<{\langle}
\def\>{\rangle}
\def\div{{\rm div}}

\def\n{\nabla}

\def\RR{\mathbb{R}}

\def\BB{\mathbb{B}}
\def\SS{\mathbb{S}}

\def\o{\omega}
\def\O{\widehat\Sigma}
\def\L{\mathcal{L}}
\def\p{\partial}
\def\e{\epsilon}
\def\a{\alpha}
\def\b{\beta}
\def\g{\gamma}

\def\s{\sigma}

\def\ep{\epsilon}

\def\W{\mathcal W}
\def\V{\mathcal V}

\def\bn{\mathbf{b}_{n+1}}
\def\bt{\mathbf{b}_\theta}
\def\Bn{{\mathbb B}^{n+1}}

\def\R{\mathbb{R}}

{\left\{\begin{array}{@{}l@{}}}{\end{array}\right.}
\patchcmd{\abstract}{\scshape\abstractname}{\textbf{\abstractname}}{}{}
\makeatletter 
\def\@makefnmark{} 
\makeatother

\numberwithin{equation}{section}

\makeatletter
\@namedef{subjclassname@2020}{\textup{2020} Mathematics Subject Classification}
\makeatother

\begin{document}
\title[Alexandrov-Fenchel inequalities]{Alexandrov-Fenchel inequalities for convex hypersurfaces in the half-space with capillary boundary}
 \author{Guofang Wang, Liangjun Weng and Chao Xia}

\address{Mathematisches Institut, Albert-Ludwigs-Universit\"{a}t Freiburg, Freiburg im Breisgau, 79104, Germany}
\email{guofang.wang@math.uni-freiburg.de}

\address{School of Mathematical Sciences, Anhui University, Hefei, 230601, P. R. China}

\email{ljweng08@mail.ustc.edu.cn}

\address{School of Mathematical Sciences, Xiamen University, Xiamen, 361005, P. R. China}
\email{chaoxia@xmu.edu.cn}

 
\maketitle
\begin{abstract}
	In this paper, we first introduce quermassintegrals for capillary hypersurfaces in the half-space. Then we solve the related isoperimetric type problems for the convex capillary hypersurfaces and obtain   the corresponding Alexandrov-Fenchel inequalities. 
	In order to prove these results, we construct a new locally constrained  curvature flow and  prove that the flow converges  globally to a spherical cap.
	
\end{abstract}


\section{Introduction}

    Let $\Sigma$ be a closed embedded hypersurface in $\RR^{n+1}$ and 
$\widehat{\S}$  the domain enclosed by $\Sigma$ in $\mathbb{R}^{n+1}$. The classical isoperimetric inequality 
states
\begin{eqnarray}
\label{01}
\frac {|\Sigma|}{\omega_n}\ge  \left( \frac {|\widehat \Sigma|}{\bn} \right)^{\frac n {n+1}},
\end{eqnarray} 
with equality holds if and only if $\S$ is a sphere. Here $\bn=|\Bn|$, the volume of the unit ball  $\Bn$, and   $\omega_n=(n+1) \bn=|\SS^n|$, the area of the unit sphere
$\SS^n$.
Its natural generalization 
is the following classical Alexandrov-Fenchel inequality 

\begin{eqnarray}\label{af ineq closed}
\frac{ \V_{l}(\widehat{\S}) }{\bn}\ge \left(\frac{ \V_{k}(\widehat{\S}) }{\bn}\right)^{\frac{n+1-l}{n+1-k}},\quad 0\le   k<l  \leq n, 
\end{eqnarray}
with equality holds if and only if $\S$ is a sphere,
provided that 
$\Sigma$ is a $C^2$ convex hypersurface.
Here  $\V_{k+1}(\widehat{\S})$ is the quermassintegral of $\widehat \Sigma$ defined by
\begin{eqnarray}\label{eq00}
\V_{k+1}(\widehat{\S}):=\frac 1{n+1} \int_{\S} H_{k }\, dA,\quad 0\leq k\leq n  ,  \quad \V_0(\widehat{\S}) = {| \widehat \Sigma|} ,
\end{eqnarray}
where $H_k$  ($1\le k\le n$) is the $k$-$th$ normalized mean curvature of $\Sigma\subset \R^{n+1}$ and $H_0=1$. It was proved in \cite{GL09} that \eqref{af ineq closed} holds true if 
$\Sigma$ is $l$-convex and star-shaped. Here by $l$-convex  we mean that $H_j>0$ for all $j\le l$.  {The case $k=0$ was proved to be true for $(l+1)$-convex  hypersurfaces in \cite{Chang, Qiu}. 
The case $l=2, k=1$, in which \eqref{af ineq closed} 
is called Minkowski's inequality, was  proved to be also true for outward minimizing sets in \cite{Hui07} (see also \cite{FS} and a very recent work \cite{AFM2} 
 by using a nonlinear potential theory.)}
 It remains still open whether \eqref{af ineq closed} is true for all $l$-convex hypersurfaces, { except  the case $l=2, k=0$ that has been  proved in \cite{AFM2}.} 
 
In this paper we are interested in its generalization to hypersurfaces with boundary. More precisely, we consider 
 hypersurfaces in $ \overline \RR^{n+1}_+$ with boundary supported on the  hyperplane $ \p \overline \RR^{n+1}_+$. Let $\S$ be a compact manifold with boundary $\p \Sigma$, which  is  properly  embedded hypersurface into  $\overline \RR^{n+1}_+$. In particular,
 $\text{int} (\Sigma)\subset \RR^{n+1}_+$ and $\p \Sigma \subset  \p\overline \RR^{n+1}_+$.
 Let $\widehat \Sigma$ be the bounded  domain enclosed by $\Sigma$ and the hyperplane $\p \overline \RR^{n+1}_+$. It is clear that the following relative isoperimetric inequality follows from \eqref{01}
\begin{eqnarray}
\label{03}
\frac {|\Sigma|}{|{\mathbb S}_+^{n}|}\ge  \left( \frac {|\widehat \Sigma|}{|{\mathbb B}^{n+1}_+|} \right)^{\frac n {n+1}},
\end{eqnarray} 
where $ {\mathbb B}^{n+1}_+$ is the upper half unit ball and ${\mathbb S}_+^{n}$ is  the 
upper half unit sphere. For the  relative isoperimetric inequality outside a convex domain, see \cite{CGR, LWW, FM}. As a domain in $\RR^{n+1}$, $\widehat \Sigma$ has  a boundary, which consists of two parts: one is $\Sigma$ and the other, which
will be denoted by $\widehat{\p \Sigma}$,  lies on   $\p \overline \RR^{n+1}_+$.
Both have a common boundary, namely $\p \Sigma$. (See Figure 1 below.) 
 Instead of just considering the area of $\Sigma$, it is interesting to consider 
the following free energy  functional
\begin{eqnarray}
\label{04}
 |\Sigma|-\cos \theta |\widehat {\p\Sigma}|
\end{eqnarray} 
for a fixed angle constant $\theta \in (0,\pi)$. The second term $\cos \theta |\widehat {\p\Sigma}|$ is the so-called wetting energy in the theory of capillarity (See for example \cite{Finn}).
 If we consider to minimize this functional under the constraint that  the volume $|\widehat \S|$ is fixed, then we have the following optimal inequality, which is called the capillary isoperimetric inequality,
\begin{eqnarray}
\label{05}
\frac {|\Sigma|-\cos \theta |\widehat{\p \Sigma}|}
{|{\mathbb S}_{\theta}^{n}|-\cos \theta |  { \widehat{\p \SS^n_\theta}}|}\ge  \left( \frac {|\widehat \Sigma|}{|{\mathbb B}^{n+1}_\theta|} \right)^{\frac n {n+1}},
\end{eqnarray} 
with equality holding if and only if {$\Sigma$ is homothetic to $\SS^n_\theta$, namely  a spherical cap \eqref{sph-cap} with contact angle $\theta$. }Here $\SS_\theta$, $\widehat {\p \SS_\theta}$, $ {\mathbb B}^{n+1}_\theta$ are defined by
$$ \SS^n_\theta =\{ x\in \SS^n \,|\,  \langle x, e_{n+1} \rangle >\cos \theta\}, \quad {\mathbb B}^{n+1}_\theta =\{ x\in   {\mathbb B}^{n+1} \,|\,  \langle x, e_{n+1} \rangle >\cos \theta\},
$$
$$ { \widehat{\p \SS^n_\theta}= \{ x\in   {\mathbb B}^{n+1} \,|\,  \langle x, e_{n+1} \rangle =\cos \theta\},
}
$$
where $e_{n+1}$ the $(n+1)$-th standard basis in $\R^{n+1}_+$. For simplicity we denote 
\begin{eqnarray} \label{06}
\textbf{b}_\theta : =  \textbf{b}^\theta_{n+1}:=|{\mathbb B}^{n+1}_\theta|,  \qquad\quad 
\omega_\theta:= \omega_{n, \theta}:=  |{\mathbb S}_{\theta}^{n}|-\cos \theta | {\widehat{ \p \SS^n_\theta}}|.
\end{eqnarray}
The explicit formulas for $\textbf{b}_\theta$ and $\omega_\theta$ will be given in Section 2.3 below.
In particular, it is easy to check that
$ (n+1) \textbf{b}_\theta =\omega_\theta.$
The proof of \eqref{05} is not trivial, which uses the spherical symmetrization. See for example \cite[Chapter 19]{Maggi}. 
For related physical problems, one can refer to the classical book of Finn  \cite{Finn}.

The main objectives of this paper are  considering the following problems
\begin{enumerate}
	\item  To find suitable generalizations of the  quermassintergals $\V_k$ for hypersurfaces with boundary supported on $\p \overline \RR^{n+1}_+$,
	which are closely related to the free energy \eqref{04}.
	\item To establish the Alexandrov-Fenchel inequality for these new
	 quermassintegrals.
\end{enumerate}

To answer  the first question, we introduce the following new geometric functionals.

\begin{eqnarray*} 
	\V_{0,\theta}(\widehat{\Sigma}):= |\widehat{ \S}|,
\qquad \V_{1,\theta}(\widehat{\Sigma}):= \frac{1}{n+1} ( |\S|-\cos\theta |\widehat{\p\S}|),\end{eqnarray*}
and for $1\leq k\leq n$,
\begin{eqnarray}\label{quermassintegrals}
\V_{k+1,\theta}(\widehat{\S}) 
&:=&
\frac{1}{n+1}\left(\int_\S H_kdA - \frac{\cos\theta \sin^k\theta  }{n}\int_{\p\S} H_{k-1}^{\p\S}ds \right)
 ,
\end{eqnarray}where  $H_{k-1}^{\p\S}$ is the normalized  $(k-1)$-th mean curvature of $\p\S\subset \RR^n$ (see Section 2 for details).
In particular, one has
\begin{eqnarray*}
	\V_{n+1,\theta}(\widehat{\S}) =\frac{1}{n+1} \int_\S H_ndA-\cos\theta\sin^n\theta \frac{\o_{n-1}}{n(n+1)} .
\end{eqnarray*}
 We believe that these quantities are the suitable  quermassintegrals for hypersurfaces with boundary intersecting with $\p \overline{ \mathbb{R}}^{n+1}_+$ at angle $\theta\in (0,\pi)$, which is supported by the following result.

\begin{thm}\label{thm0}
	Let $\S_t\subset \overline{{\R}}_+^{n+1}$ be a family of smooth, embedded capillary  hypersurfaces with a constant contact angle $\theta\in (0,\pi)$, which are  given by the embedding $x(\cdot,t):M\to \overline{\R}^{n+1}_+$ and satisfy
	\begin{eqnarray*}
		(\p_t x)^\perp =f\nu,
	\end{eqnarray*}for some speed function $f$. Then for $0\leq k\leq n$,
	\begin{eqnarray}\label{variation}
		\frac{d}{dt} \V_{k,\theta}(\widehat{\S_t})=\frac{n+1-k}{n+1}\int_{\S_t} fH_{k}dA_t,
	\end{eqnarray}
	and in particular \begin{eqnarray} \label{GB00}
		\frac{d}{dt}\V_{n+1,\theta}(\widehat{\S_t})=0.
	\end{eqnarray}
\end{thm}

 A hypersurface in $\overline{{\R}}_+^{n+1}$  with boundary supported on $\p\overline{{\R}}_+^{n+1}$ is called capillary hypersurface if
it intersects with $\p\overline{{\R}}_+^{n+1}$ at a constant angle.
 For closed hypersurfaces in $\R^{n+1}$,  a similar variational formula
as \eqref{variation} characterizes the quermassintegrals in \eqref{eq00}. This formula  for the quermassintegrals is also true for closed hypersurfaces in other space forms. See for example \cite{WX14}.

Our second result is the generalized Alexandrov-Fenchel inequalities for convex capillary hypersurfaces. 

\begin{thm}\label{thm 1}
	For $n\ge 2$,  let $\Sigma\subset    \overline{\mathbb{R}}^{n+1}_+$ be a convex capillary hypersurface with a constant contact angle $\theta \in (0, \frac{{\pi}}{2} ]$, then 
	  there holds
	\begin{eqnarray}\label{af ineq}
	\frac{ \V_{n,\theta}(\widehat{\S}) }{\textbf{b}_\theta}
	\geq 
	\left( \frac{  \V_{k,\theta}(\widehat{\S}) }{\textbf{b}_\theta} \right)^{\frac{1}{n+1-k}} ,  \quad  \forall \,  0\le k <n,
	\end{eqnarray}
	with equality if and only if $\Sigma$  is a spherical cap in \eqref{sph-cap}. 	Moreover,
		\begin{eqnarray} \label{GB}
	\V_{n+1,\theta}(\widehat{\S})=\omega_\theta=(n+1)\textbf{b}_\theta,
	\end{eqnarray} 	
\end{thm}

\eqref{GB} follows easily from \eqref{GB00}, by constructing a smooth family of capillary hypersurfaces connecting to a spherical cap given in \eqref{sph-cap}. Therefore \eqref{GB} 
is true for any  capillary hypersurfaces with contact angle $\theta$.
Moreover, it is equivalent to
\begin{eqnarray}\label{GB10}
	 \int_\S H_ndA=(n+1)\omega_\theta +\cos\theta\sin^n\theta \frac{\o_{n-1}}{n} =  {|\SS^n_\theta|} ,
\end{eqnarray}
a Gauss-Bonnet type result for capillary hypersurfaces with contact angle $\theta$. When $n=2$, \eqref{GB10} implies a Willmore   inequality for  capillary hypersurfaces with contact angle $\theta$.
\begin{cor}[Willmore]\label{cor 1.3}
	Let   $\Sigma\subset \overline{\mathbb{R}}^3_+$ be a convex   capillary surface with a constant contact angle $\theta \in (0,\frac{\pi}{2}]$, then 
\begin{eqnarray}\label{Willmore}
\int_\Sigma  H^2 dA \ge 4 |\SS^2_\theta|,
\end{eqnarray}with equality holds if and only if $\S$ is spherical cap in \eqref{sph-cap}.
\end{cor}
Here $H=2H_1$ is the ordinary mean curvature for surfaces and it is obvious that $H^2 \ge 4H_2$.  
When  $n=2$, \eqref{af ineq} implies  a Minkowski  type inequality for convex capillary  surfaces with    boundary in $\overline{\mathbb{R}}^3_+$.
\begin{cor}[Minkowski]\label{cor 1}	Let $\Sigma\subset \overline{\mathbb{R}}^3_+$ be a convex  capillary surface with a constant contact angle  $\theta \in (0,\frac{\pi}{2}]$, then
	\begin{eqnarray}\label{minkowski ineq}
	\int_\S H dA  \geq 2 \sqrt{\omega_{2,\theta}} \cdot (|\S|-\cos\theta |\widehat{\p\S}|)^{\frac{1}{2}}+ \sin\theta \cos\theta |\p\S|,
	\end{eqnarray} 
where	\begin{eqnarray*}
	\omega_{2,\theta}= 3	\bt= \left(2-3\cos\theta +\cos^3\theta\right)\pi.
	\end{eqnarray*} 
	Moreover,  equality holds if and only if $\Sigma$ is a spherical cap in \eqref{sph-cap}. 
\end{cor}

From these results, it is natural to propose

\

\noindent{\bf Conjecture 1.5.}   {\it 
For $n\ge 2$,  let $\Sigma\subset    \overline{\mathbb{R}}^{n+1}_+$ be a convex capillary hypersurface with a contact angle $\theta \in (0, {{\pi}})$, there holds
	\begin{eqnarray}\label{af ineq2}
	\frac{ \V_{l,\theta}(\widehat{\S}) }{\textbf{b}_\theta}
	\geq 
	\left( \frac{  \V_{k,\theta}(\widehat{\S}) }{\textbf{b}_\theta} \right)^{\frac{n+1-l}{n+1-k}} ,  \quad  \forall 0\le k <l\le n,
	\end{eqnarray}
	with equality iff $\Sigma$ is  a spherical cap in \eqref{sph-cap}.
}

\

It would be also interesting to ask further if the conjecture is true for $k$-convex capillary hypersurfaces. 
In order to prove the conjecture, we introduce a suitable nonlinear curvature flow, which preserves $\V_{l,\theta}(\widehat{\S})$ and increases $\V_{k,\theta}(\widehat{\S})$, see Section \ref{sect 3}. If the flow globally converges to a spherical cap, then
we have  the general Alexandrov-Fenchel inequaltiy \eqref{af ineq2}. However, due to technical  difficulties we are only able to prove  in this paper the global convergence, for $l=n$ and $\theta\in (0, \frac{\pi}{2}]$, namely, Theorem \ref{thm 1}.

To be more precise, let us first recall the related work on the proof of Alensandrov-Fenchel inequalities by using  geometric curvature flow for closed hypersurfaces  in $\R^{n+1}$.
When $\p\S=\emptyset$, and we denote $\widehat{\S}$ be the bounded convex domain enclosed by $\Sigma$ in $\mathbb{R}^{n+1}$. In convex geometry,
the Alexandrov-Fenchel inequalities \eqref{af ineq closed}  between quermassintegrals  $\V_k$ and $\V_l$ play  an important role. 
In fact there are more general inequalities. See \cite{Al1, Al2, Sch} for instance. 
It is an interesting question, if
one can use a curvature flow to reprove such inequalities.
In \cite{Mc}, McCoy introduced a normalized nonlinear curvature flow to reprove the Alexandrov-Fenchel inequalities \eqref{af ineq closed} for convex domains in Euclidean space. Later, Guan-Li \cite{GL09}  weakened  the convexity condition and  only assumed that the closed hypersurface $\S$ is 
$k$-convex and star-shaped by using the inverse curvature flow, which is defined by 
\begin{eqnarray}
	\p_t x=\frac{H_{k-1}}{H_k}\nu.
\end{eqnarray}
This flow was previously studied by Gerhardt \cite{Gerhardt} and Urbas \cite{Urb90}. One of key observations in the study of this flow is that the $k$-convexity and star-shaped are preserved along this flow. This flow is also equivalent to the rescaled one \begin{eqnarray}\label{rescaled guan-li flow}
	\p_t x=\left(\frac{H_{k-1}}{H_k}-\langle x,\nu\rangle \right)\nu,
\end{eqnarray}see \cite{GL18,GL21} for instance. The motivation to use such a flow \eqref{rescaled guan-li flow} is its nice properties that
  the quermassintegrals $\V_k(\widehat{\S})$ is preserved and $\V_{k+1}(\widehat{\S})$ is non-decreasing along this flow, which follows from the  well-known Minkowski 
  formulas. Using similar geometric flows, there have been a lot of work to establish
  new Alexandrov-Fenchel inequalities  in the hyperbolic space \cite{AW18,AHL,BGL,DF,GWW,HL,HLW,LWX,SX,WX14} and  in the sphere \cite{CGLS,ChS,MS,WeiX}. 
 For the anisotropic analogue of Alexandrov-Fenchel (Minkowski) type inequalities we refer to \cite{AW21,WeiX21,Xia}.

If $\p\S\neq \emptyset$, the study of geometric inequalities with free boundary or general capillary boundary has attracted much attention in the last decades.
 For  related  relative isoperimetric inequalities and the Alexandrov-Fenchel inequalities, see for instance \cite{BoSp,BuMa,CRS,FI13,LS17,SWX,WeX21} etc. 
 Recently in \cite{SWX} Scheuer-Wang-Xia  introduced the definition of quermassintegrals  for hypersurfaces with free boundary in the Euclidean unit ball $\overline{\BB}^{n+1}$ from the viewpoint of the first variational formula, and they proved the highest order Alexandrov-Fenchel inequalities  for convex hypersurfaces with free boundary  in $\overline{\BB}^{n+1}$. Very recently, the second and the third authors  \cite{WeX21} generalized the work in \cite{SWX} by 
  introducing  the corresponding  quermassintegrals  for   general capillary hypersurfaces and established   Alexandrov-Fenchel inequalities for convex  capillary hypersurfaces  in $\overline{\BB}^{n+1}$. The flows introduced to establish these inequalities  in \cite{SWX, WeX21} are motivated by  new Minkowski formulas proved in \cite{WX1}.
  
 Now we introduce our curvature flow for capillary hypersurfaces in the half space. Let $e=-e_{n+1}$, where $e_{n+1}$ the $(n+1)$-th coordinate in $\R^{n+1}_+$. 
 Let $x:\Sigma \to  \overline{\mathbb{R}}^{n+1}_+$ with boundary $x_{|\p\Sigma}:\partial \Sigma \to \partial \RR_+^{n+1}$ and $\nu$ its unit normal vector field. We introduce
 \begin{eqnarray}\label{flow0}
\left( \p_t x \right)^\perp = \left[ \left(1+\cos \theta \<\nu, e\> \right)\frac {H_{l-1}}{H_l} -\<x, \nu\>\right] \nu.
 \end{eqnarray}
 Using the Minkowski formulas given in \eqref{minkowski sigma_k},  we show that flow \eqref{flow0} preserves $\V_{l,\theta}(\widehat{\S})$, while increases $\V_{k,\theta}(\widehat{\S})$ for $k<l$. However, due to the weighted function  in the flow we are only able at moment to show that flow \eqref{flow0} preserves the convexity, when $l=n$. In this case, we can further bound $\frac{H_n}{H_{n-1}}$. In order to 
 bound all principal curvature we need to estimate the mean curvature, which satisfies a nice evolution equation \eqref{evo_H}. However the normal derivative of $H$, $\n_{\mu}H$, has a bad sign, if $\theta > \frac{\pi}{2}$. Hence we have to restrict ourself on the range $\theta \in (0,\frac{\pi}{2}]$. Under these conditions we then succeed to show the global convergence, and hence Alexandrov-Fenchel inequalities.
 It would be interesting to ask if one can also prove the global convergence for the case $\theta> \frac {\pi}2$. In analysis, this case is related to the worse case in the Robin boundary problem for the corresponding  PDE.
  We remark also that  a  convex capillary hypersurface with contact angle $\theta\leq  \frac{\pi}{2}$ could have different geometry from that with  $\theta > \frac{\pi}{2}$. The former was called a convex cap and was studied in \cite{Busemann}.

Comparing with inequalities  established in \cite{SWX, WeX21}, which  are actually  implicit inequalities and  involve inverse functions of certain geometric quantities  that can not be explicitly expressed by elementary functions, we have here a  geometric inequality \eqref{af ineq} in an explicit and clean form. An optimal inequality with an  explicit form has more applications. 
A further  good example was given very recently in an  optimal insulation problem in \cite{DNT}, where the optimal inequalities between $\V_n$ and $\V_k$ for any $k<n$ for closed hypersurfaces have been used crucially. We expect that our results can be similarly used 
in an  optimal insulation problem  for capillary hypersurfaces.

    \
 
 \textit{The rest of the article is structured as follows.} In Section \ref{sect 2}, we introduce the   quermassintegrals for capillary  hypersurfaces and collect the relevant evolution equations to finish the proof of Theorem \ref{thm0}. In Section \ref{sect 3}, we introduce our nonlinear inverse curvature flow and show the monotonicty of our quermassintegrals \eqref{quermassintegrals} under the flow. In Section \ref{sect 4}, we obtain  uniform estimates for convex capillary  hypersurfaces along the flow 
  and the global convergence.  Section \ref{sect 5} is devoted to prove the Alexandrov-Fenchel inequalities  for convex  capillary hypersurfaces  in the half-space, i.e. Theorem \ref{thm 1}.
\section{Quermassintegrals and Minkowski formulas}\label{sect 2}

 Since we will deform hypersurfaces by studying a geometric flow, it is convenient to use immersions.
Let $M$ denote a compact orientable smooth manifold of dimension $n$ with  boundary $\p M$, and $x:M\to \overline{\mathbb{R}}^{n+1}_+$ be a proper smooth immersed hypersurface. In particular, $x(\text{int}(M))\subset \RR^{n+1}_+$ and $x(\p M)\subset \p\RR^{n+1}_+$. 
Let $\S=x(M)$ and $\partial \S=x(\partial M)$. If no confusion, we will do not distinguish the hypersurfaces $\S$ and the immersion $x:M \to \overline{\R}^{n+1}_+$. {Let $\O$ be the  bounded domain enclosed by $\S$ and $\p{\RR}^{n+1}_+$.
Let $\nu$ and $\overline{N}$ be the unit outward normal  of $\S\subset {\widehat\Sigma}$ and $\p{\RR}^{n+1}_+ \subset \overline {\RR}^{n+1}_+$ respectively.}
\subsection{Higher order mean curvatures}\label{sec 2.1}

For $\k=(\k_1,\k_2\cdots, \k_n) \in \RR^n$, let  $\s_k(\k), k=1,\cdots, n,$ be the $k$-th elementary symmetric polynomial functions and $H_k(\k)$ be its normalization $H_k(\k)=\frac{1}{\binom{n}{k}}\s_k(\k)$.  For $i=1,2, \cdots, n$,  let $\k|i\in \R^{n-1}$ 
(or $\kappa | \kappa_i$) denote $(n-1)$ tuple deleting the $i$-th component from $\k$.

We shall use the following basic properties about $\sigma_k$.
\begin{prop}\label{prop2.2}\
	\begin{enumerate}
		\item $ \sigma_k(\k)=\sigma_k(\k|i)+\k_i\sigma_{k-1}(\k|i),  \quad  \forall 1\leq i \leq n.$
		\item $ \sum\limits_{i = 1}^n {\sigma_{k}(\k|i)}=(n-k)\sigma_k(\k).$
		\item $\sum\limits_{i = 1}^n {\k_i \sigma_{k-1}(\k|i)}=k\sigma_k(\k)$.
		\item $\sum\limits_{i = 1}^n {\k_i^2 \sigma_{k-1}(\k|i)}=\s_1(\k)\sigma_k(\k)-(k+1)\s_{k+1}(\k)$.
	\end{enumerate}
\end{prop}
Let  $ \Gamma_+:=\{\kappa\in\mathbb{R}^n:\kappa_i>0, 1\leq i\leq n\}$ and $\Gamma_+^k=\{\kappa\in \R^n \,| \, H_j(\kappa)>0, \forall 1\leq j\le k\}$. It is clear that $\Gamma_+=\Gamma_+^n$.
\begin{prop}\label{prop2.3}
	For $1\le k<l\le n$,  we have
	\begin{eqnarray}\label{NM-ineq}
		H_kH_{l-1}\geq H_{k-1}H_{l }, \quad \quad \forall\quad  \kappa\in  \Gamma^l_+,
	\end{eqnarray}	
	with equality holds if and only if $\kappa=\lambda(1, \cdots, 1)$ for any $\lambda>0$.
	Moreover, $F(\k)=\frac{\s_k}{\s_{k-1}}(\k)$ is concave in $\Gamma^k_+$.
\end{prop}

These are well-known properties. For a proof we refer to \cite[Chapter XV, Section 4]{Lie} and \cite[Lemma 2.10, Theorem 2.11]{Spruck} respectively.

We use $D$ to denote the Levi-Civita connection of $\overline{\RR}^{n+1}_+$ with respect to the Euclidean metric $\delta$, 
and $\n$  the Levi-Civita connection on $\S$ with respect to the induced metric $g$ from the immersion $x$. 
The operator $\div, \Delta$, and $\n^2$ are  the divergence, Laplacian, and Hessian operator on $\S$ respectively.
The second fundamental form $h$ of $x$  is defined by
$$D_{X}Y=\n_{X}Y- h(X,Y)\nu.$$ 

Let $\k=(\k_1, \k_2,\cdots, \k_n)$ be the set of principal curvatures, i.e, the set of eigenvalues of $h$. Then we
 denote  $\sigma_k=\sigma_k(\k)$ and $H_k=H_k(\k)$ resp. be  the $k$-th mean curvature and 
the normalized $k$-th mean curvature of $\Sigma$.
 We  also use the convention that $$\s_0=H_0=1 \quad \s_{n+1}=H_{n+1}=0.$$

\begin{rem}\label{tensor}
	We will simplify the notation by using the following shortcuts occasionally:
	\begin{enumerate}
		\item When dealing with complicated evolution equations of tensors, we will use a local frame to express tensors with the help of their components, i.e. for a tensor field $T\in \mathcal{T}^{k,l}(\S)$, the expression $T^{i_{1}\dots i_{k}}_{j_{1}\dots j_{l}}$ denotes  
		\begin{eqnarray*}
			T^{i_{1}\dots i_{k}}_{j_{1}\dots j_{l}}=T(e_{j_{1}},\dots,e_{j_{l}},\e^{i_{1}},\dots \e^{i_{k}}),
		\end{eqnarray*}
		where $(e_{i})$ is a local frame and $(\e^{i})$ its dual coframe. 
		\item { The $m$-th covariant derivate  of  a $(k,l)$-tensor field $T$, $\nabla^m T$,  is locally expressed by 
		\begin{eqnarray*}
			T^{i_1\dots i_k}_{j_1\dots j_l; j_{l+1}\dots j_{l+m}}.
		\end{eqnarray*}
}	
		\item We shall use the convention of the Einstein summation. For convenience the components of the Weingarten map $\W$ are denoted by $(h^{i}_{j})=(g^{ik}h_{kj})$, and $|h|^2$ be the norm square of the second fundamental form, that is $|h|^2=g^{ik}h_{kl}h_{ij}g^{jl}$, where  $(g^{ij})$ is the inverse of $(g_{ij})$. We use the metric tensor $(g_{ij})$ and its inverse $(g^{ij})$ to  lower down and raise up the indices of tensor fields on $\S$.
		
	\end{enumerate}
\end{rem}

\subsection{Quermassintegrals in the half-space}
In order to introduce our  quermassintegrals for capillary  hypersurfaces in the half-space, we review first the quermassintegrals in $\RR^{n+1}$, see e.g. \cite{Sch}. Given a bounded convex domain $\O\subset \mathbb{R}^{n+1}$ with smooth boundary $\p\O$, its $k$-th  quermassintegral is defined by
\begin{eqnarray*}
	\V_0(\O):=|\O|,
\end{eqnarray*}
and for $0\leq k\leq n$,
\begin{eqnarray*}
	\V_{k+1}(\O):=\frac{1}{n+1} \int_{\p \O}H_{k}  dA,
\end{eqnarray*}
where $H_{k} $ is the normalized $k$-th mean curvature of $\p\O\subset \R^{n+1}$.  One can check that
\begin{eqnarray}\label{add_1}
\frac d {dt} \V_{k+1}(\O_t) =\frac{n-k}{n+1} \int_{\partial \widehat\Sigma_t} H_{k+1} f dA,
\end{eqnarray}
for  a family of bounded convex bodies $\{\O_t\}$ in $\R^{n+1}$ whose boundary $\p\O_t$ evolving by a  normal variation  with speed function $f$.  For a proof see e.g. \cite[Lemma 5]{GL09}. As mentioned above, a similar first variational formula also holds
in space forms, see \cite{Reilly}. Therefore formula \eqref{add_1}
is the characterization of the quermassintegrals for closed hypersurfaces in space forms.

Now we define the following geometric functionals for convex hypersurfaces $\S$ with capillary boundary in $\overline{\mathbb{R}}^{n+1}_+$ with a constant contact angle $\theta$ along $\p\S\subset \R^n$.
Let
\begin{eqnarray*} 
	\V_{0,\theta}(\widehat{\Sigma}):= |\widehat{ \S}|,
\end{eqnarray*}
 \begin{eqnarray*} 	\V_{1,\theta}(\widehat{\Sigma}):= \frac{1}{n+1} \left( |\S|-\cos\theta |\widehat{\p\S}|\right).\end{eqnarray*}
and for $1\leq k\leq n$,
\begin{eqnarray*} 
	\V_{k+1,\theta}(\widehat{\S}): = \frac{1}{n+1}\left(\int_\S H_kdA -\frac{\cos\theta \sin^k\theta  }{n}\int_{\p\S} H_{k-1}^{\p\S}ds \right).
\end{eqnarray*}{Here $H_{k-1}^{\p\S}:=\frac{1}{\binom{n-1}{k-1}}\s_{k-1}(\hat{\k})$ is the normalized  $(k-1)$-th mean curvature of $\p\S\subset \RR^n$ and $\s_{k-1}(\hat{\k})$ is the $(k-1)$-elementary symmetric function on $\RR^{n-1}$ evaluating at the principal curvatures $\hat{\k}$ of $\p\S\subset \RR^n$. } In particular, we have
\begin{eqnarray*}
	\V_{2,\theta}(\widehat{\S})=\frac{1}{n(n+1)}\left(\int_\S H dA-\sin\theta \cos\theta | {\p\S}|\right).
\end{eqnarray*} Here $H$ is the (un-normalized) mean curvature, i.e. $H=n H_1$.
From Gauss-Bonnet-Chern's theorem, we know
\begin{eqnarray*}
	\V_{n}(\O)=\frac{\omega_{n-1}}{n},
\end{eqnarray*}if $\O\subset \RR^{n}$ is a convex body (non-empty, compact, convex set).  As a result, we see
\begin{eqnarray*}
	\V_{n+1,\theta}(\widehat{\S}) =\frac{1}{n+1} \int_\S H_ndA-\cos\theta\sin^n\theta \frac{\o_{n-1}}{n(n+1)} .
\end{eqnarray*}

\subsection{Spherical caps}
Let $e:=-e_{n+1}=(0,\ldots,0,-1)$.
We consider	a family of spherical caps lying entirely in   $\overline{\mathbb{R}}^{n+1}_+$ and intersecting $\mathbb{R}^{n}$ with a constant contact angle $\theta\in (0,\pi)$ given by
	\begin{eqnarray}\label{sph-cap}
  C_{r,\theta }(e) :=\Big\{x\in \overline{\mathbb{R}}^{n+1}_+ \big| |x-r\cos\theta e|= r\Big\},\quad r\in [0,\infty),
	\end{eqnarray} which has radius $r$ and centered at   $ r\cos\theta e $. {To emphasize $e$ and to distinguish with the center of the spherical cap, $r\cos\theta e$,  we call $ C_{r,\theta }(e)$ 
	a spherical cap around $e$.}		
If without confusion, we just write $C_{r,\theta}$ for $C_{r,\theta}(e)$ in the rest of this paper.
One can easily check that $  C_{r,\theta}$ is the static solution to flow \eqref{flow with capillary} below, that is,
\begin{eqnarray}\label{static model radius}
	1+  \cos\theta \langle \nu,e\rangle -\frac{1}{r}\langle x,\nu\rangle=0,
\end{eqnarray}
and it intersects with the support $\partial \overline{\mathbb{R}}^{n+1}_+$ at the  constant angle $\theta$.%

The volume of $\widehat {C_{r,\theta}}$ 
$$\V_{0,\theta} (\widehat {C_{r,\theta}})= r^{n+1} \bt,$$
where $\bt$ is the volume of $C_{1,\theta}$, which is congruent to $\SS_\theta$, defined in the introduction. 
One can compute 
 \begin{eqnarray}
 \label{bt}
 \bt:=\frac{\o_n}{2} I_{\sin^2\theta}(\frac{n}{2},\frac{1}{2})-\frac{\o_{n-1}}{n} \cos\theta \sin^n\theta,
  \end{eqnarray} and $I_s(\frac{n}{2},\frac{1}{2})$ is the regularized incomplete beta function given by 
	\begin{eqnarray}
	I_s(\frac{n}{2},\frac{1}{2}):=\frac{\int_0^s t^{\frac{n}{2}-1} (1-t)^{-\frac{1}{2} }dt}{\int_0^1 t^{\frac{n}{2}-1} (1-t)^{-\frac{1}{2} }dt}.
	\end{eqnarray} 
	Moreover, one can readily  check that
$$\V_{1,\theta} (\widehat {C_{r,\theta}})=\frac 1{n+1} (|C_{r,\theta}|-\cos \theta |\widehat{\partial C_{r,\theta}}|)=r^n \bt$$
and
$$\V_{k,\theta} (\widehat {C_{r,\theta}}) = r^{n+1-k} \bt.$$
Therefore, $C_{r,\theta}$ achieves  equality in the Alexandrov-Fenchel inequalities \eqref{af ineq}.

\subsection{Minkowski formulas} As above, $\S\subset \overline{\RR}^{n+1}_+$ is  a smooth, properly embedded  capillary hypersurface, given by the embedding $x:M\to \overline{\RR}^{n+1}_+$, where $M$ is a compact, orientable smooth manifold of dimension $n$ with non-empty boundary. Let $\mu$ be the unit outward co-normal of $\p\S$ in $\S$ and  $\overline{\nu}$ be the unit normal to $\partial\Sigma$ in $\partial\mathbb{R}^{n+1}_+$ such that $\{\nu,\mu\}$ and $\{\overline{\nu},\overline{N}\}$ have the same orientation in normal bundle of $\partial\Sigma\subset\overline{\mathbb{R}}^{n+1}_+$. We define 
the contact angle $\theta$ between  the hypersurface $\Sigma$ and the support $\partial\overline{\R}^{n+1}_+$ by
$$\< \nu, \overline N\>=\cos (\pi-\theta).$$
 It  follows  
\begin{eqnarray}\label{co-normal bundle}
	\begin{array} {rcl}
		\overline{N} &=&\sin\theta \mu-\cos\theta \nu,
		\\
		\overline{\nu} &=&\cos\theta \mu+\sin\theta \nu,
	\end{array}
\end{eqnarray}
or equivalently
\begin{eqnarray}\label{2.5}
\begin{array} {rcl}
\mu &=&\sin\theta\overline N + \cos \theta \overline \nu ,
\\
{\nu} &=&-\cos \theta \overline N +\sin\theta \overline \nu.
\end{array}
\end{eqnarray}
 $\p\S$ can be viewed as a smooth closed hypersurface in $\RR^n$, which bounds a bounded domain $\widehat{\p\S}$ inside $\RR^n$. By our convention, $\bar \nu$ is the unit outward normal of $\p \S$ in  $\widehat{\p\S}\subset \RR^n$.  See Figure 1.
  \begin{center}
  	\begin{figure}[h] \includegraphics[width=0.85\linewidth]{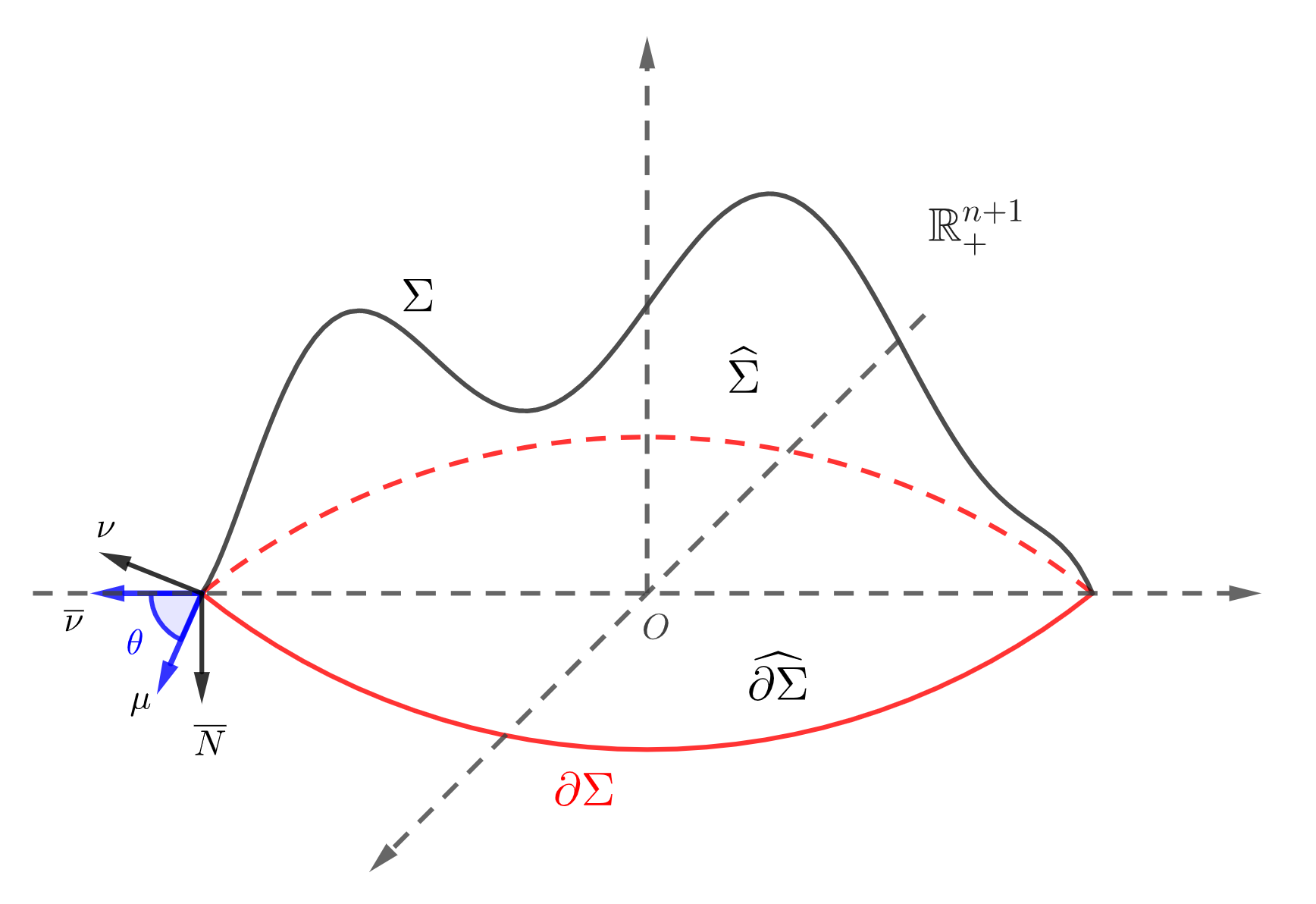} \label{fig1}  \caption*{Figure 1. A capillary hypersurface $\S$ with a contact angle $\theta$}   
  	\end{figure}
  \end{center}

 The second fundamental form of $\p \S$ in  $\RR^n$ is given by $$\widehat{h}(X, Y):= -\<\nabla^{\RR^n}_X Y, \bar \nu\>= -\<D_X Y, \bar \nu\>, \quad X, Y\in T(\p\S).$$
The second equality holds since $\<\bar \nu, \bar N\circ x\>=0$.
The second fundamental form of $\p \S$ in $\S$ is given by
$$\widetilde{h}(X, Y):= -\<\nabla_X Y, \mu\>= -\<D_X Y, \mu\>, \quad X, Y\in T(\p\S).$$
The second equality holds since $\<\nu, \mu\>=0$.

\begin{prop}\label{basic-capillary} Let $\S\subset \overline{\RR}^{n+1}_+$ be a capillary hypersurface. Let $\{e_\a\}_{\a=2}^{n }$ be an orthonormal frame of $\p \S$. Then along $\p\S$,
	\begin{itemize}
		\item[(1)] $\mu$ is a principal direction of $\S$, that is, $h_{\mu \a}=h(\mu, e_\a)=0$.
		\item[(2)] $h_{\a\b}=\sin\theta \widehat{h}_{\a\b} .$
		\item[(3)] $\widetilde{h}_{\a\b}=\cos\theta \widehat{h}_{\a\b} =  \cot\theta h_{\a\b}.$
		\item[(4)] $h_{\a\b; \mu}=\widetilde{h}_{\b\g}(h_{\mu\mu}\delta_{\a\g}-h_{\a\g})$.
	\end{itemize}
\end{prop}

\begin{proof}The first assertion is well-known, 
	see e.g.  \cite{RR}. 
(2) and  (3) follow from 
	\begin{eqnarray*}
		h_{\alpha\beta}=-\langle D_{e_\alpha}e_\beta,\nu\rangle	=\langle \widehat{h}_{\alpha\beta}\overline{\nu} ,\nu\rangle
		=\sin\theta \widehat{h}_{\alpha\beta} ,
	\end{eqnarray*}
and
	\begin{eqnarray*}
		\widetilde{h}_{\a\b}= -\<D_{e_\alpha}e_\beta, \mu\>= \<\widehat{h}_{\alpha\beta}\overline{\nu} ,\mu\>=\cos\theta \widehat{h}_{\a\b} .
	\end{eqnarray*}
	For (4),  taking derivative of $h(\mu, e_\a)=0$ with respect to $e_\b$
	and  using the Codazzi equation and (1),  we get
	\begin{eqnarray*}
		0=e_\b(h(\mu, e_\a))&=&h_{\a \mu; \b}+h(\n_{e_\b}e_\a, \mu)+h(\n_{e_\b}\mu, e_\a)
		\\&=&h_{\a \b; \mu}+ \<\n_{e_\b}e_\a, \mu\>h_{\mu\mu}+\<\n_{e_\b}\mu, e_\g \>h_{\a\g}
		\\&=&h_{\a \b; \mu}-\widetilde{h}_{\b\g}(h_{\mu\mu}\delta_{\a\g}-h_{\a\g}).
	\end{eqnarray*}
\end{proof}
The Proposition \ref{basic-capillary} has a direct conseqeunce.
{ \begin{cor}\label{x-conv}
If $\S$ is a convex capillary hypersurface, then $\p \Sigma \subset \p \R^{n+1}$ is also convex, i.e., $\widehat h\ \ge 0$, while 
 $\p \Sigma \subset \Sigma$ is convex  $(\widetilde h \ge 0)$ if $\theta\in (0,\frac{\pi}{2}]$ and  concave $(\widetilde h \le 0)$ if $\theta \in [\frac {\pi} 2, \pi)$.
\end{cor}}

 The following Minkowski type formulas for capillary hypersurfaces  play an important role in this paper.
\begin{prop}\label{minkowski formulae}Let $x:M\to  \overline{\mathbb{R}}^{n+1}_+ $ be an smooth immersion of $\S:=x(M)$ into the half-space, whose boundary intersects $\mathbb{R}^n$ with a constant contact angle $\theta\in (0,\pi)$ along $\p\S$. 
For $1\leq k\leq n$, it holds
		\begin{eqnarray}\label{minkowski sigma_k}  
			\int_\S H_{k-1}  ( 1+\cos\theta \langle \nu,e \rangle )dA=\int_\S H_k  \langle x,\nu\rangle dA,
		\end{eqnarray}
	where $dA$ is the area element of $\S$ w.r.t.  the induced metric $g$.
\end{prop}
When $k=1$ or $\theta=\frac{\pi}{2}$, formula  \eqref{minkowski sigma_k} is known. See e.g.  \cite[Proof of Theorem 5.1]{AS} and \cite[Proposition 2.5]{GL15}. For our purpose, we need the  high order Minkowski type formulas for general $\theta$.
\begin{proof} 
Denote $  x^T :=x-\langle x,\nu\rangle \nu   $ be the tangential projection of $x$ on $\Sigma$, and $$P_e:=\langle \nu,e \rangle x-\langle x,\nu\rangle e .$$ From a direct computation, we have
	\begin{eqnarray}\label{x T}
	D_{e_i} \langle x,e_j\rangle= g_{ij}-\langle x,\nu\rangle h_{ij}
\end{eqnarray}
and		\begin{eqnarray}\label{p a T}
	\nabla_{i} (P_e^T)_{j} =\langle \nu,e \rangle g_{ij}+h_{il}\langle e,e_l\rangle \langle x,e_j\rangle-h_{il}\langle x,e_l\rangle \langle e,e_j\rangle . 
\end{eqnarray}
Along   $\partial \Sigma\subset \partial \mathbb{R}^{n+1}_+$, using \eqref{co-normal bundle}  we see  	\begin{eqnarray*}
		\langle P_e^T,\mu\rangle&=&\langle P_e,\mu\rangle =\langle \nu,e\rangle \langle x,\mu\rangle -\langle x,\nu\rangle \langle e,\mu\rangle\\&=&\langle x,-\sin\theta \nu-\cos\theta \mu\rangle,
	\end{eqnarray*}
which follows 
	\begin{eqnarray}\label{mu ortho to X P}
		\langle x^T+\cos\theta P_e^T,\mu\rangle &=&\langle x,\mu\rangle -\cos\theta\langle x,\sin\theta \nu+\cos\theta\mu\rangle \notag
		\\&=&\sin\theta \langle x,e\rangle =0. 
	\end{eqnarray}
 Denote $\sigma_{k-1}^{ij}:=\frac{\partial\sigma_k}{\partial  h^i_j }$ be the $k$-th Newton transformation.  
 Taking contraction with \eqref{x T}, \eqref{p a T} and  using Proposition \ref{prop2.2} we obtain	
	\begin{eqnarray*}
		\sigma^{ij}_{k-1}\cdot	\n_i \left((x^T+\cos\theta P_e^T)_j\right) &=&\sigma_{k-1}^{ij}
		(  g_{ij}-h_{ij}\langle x,\nu\rangle+\cos\theta\langle \nu,e\rangle g_{ij}  )
		\\&=&(n-k+1)\sigma_{k-1}( 1+\cos\theta \langle \nu,e\rangle)  -k\sigma_k  \langle x,\nu\rangle
		\\&=& \frac{n!}{(k-1)! (n-k)!} \left(
		 H_{k-1}( 1+\cos\theta \langle \nu,e\rangle)  -H_k   \langle x,\nu\rangle \right).
	\end{eqnarray*}
Using integration by parts,  we have
	\begin{eqnarray*} 
		\int_\S 	\n_i\big( \sigma^{ij}_{k-1}(x^T+\cos\theta P_e^T)_j\big)dA  = 	\int_{\p \S} \sigma_{k-1}^{ij} (x^T+\cos\theta P_e^T)_j\cdot \mu_ids .
	\end{eqnarray*}
	From   \eqref{mu ortho to X P}, we know that $(x^T+\cos\theta P_e^T)\perp \mu$ along $\p\S$.
	Since $\mu$ is a principal direction of $\S$ by Proposition  \ref{basic-capillary}, 
	we have  $  \sigma_{k-1}^{ij} (x^T+\cos\theta P_e^T)_j\cdot \mu_i=0$ along $\p\S$. It is well-known that the
	Newton tensor is divergence-free, i.e., $\n_i\sigma_{k-1}^{ij}=0$.
	Altogether yields the conclusion. 
\end{proof}

\subsection{Variational formulas}

The following first variational formula motivates us to define the quermassintegrals for capillary hypersurfaces as  \eqref{quermassintegrals}.
\begin{thm}\label{1st variation of total sigma k}
	Let $\S_t\subset \overline{{\R}}_+^{n+1}$ be a family of smooth capillary  hypersurfaces supported by  $\p\R^{n+1}_+$ with a constant contact angle $\theta\in (0,\pi)$ along $\p\S_t$, given by the embedding $x(\cdot,t):M\to \overline{\R}^{n+1}_+$, and satisfying
\begin{eqnarray}\label{normal flow}
	(\p_t x)^\perp =f\nu,
\end{eqnarray}for a smooth function $f$. Then for $-1\leq k\leq n-1$,
\begin{eqnarray}\label{fvf k+1}
	\frac{d}{dt} \V_{k+1,\theta}(\widehat{\S_t})=\frac{n-k}{n+1}\int_{\S_t} fH_{k+1}dA_t,
\end{eqnarray}
and\begin{eqnarray*}
	\frac{d}{dt}\V_{n+1,\theta}(\widehat{\S_t})=0.
\end{eqnarray*}
\end{thm}
Before proving Theorem \ref{1st variation of total sigma k}, we remark that if $\S_t\subset \overline{\RR}^{n+1}_+$ is a family of smooth capillary hypersurfaces evolving by \eqref{normal flow}, then the tangential component $(\p_tx)^T$ of $\p_t x$, which we denote by $T\in T\S_t$,  must satisfy \begin{eqnarray}\label{choiceT1}
	T |_{\p \S_t}=f\cot\theta \mu + \widetilde{T},
\end{eqnarray}where $\widetilde{T}\in T(\p\S_t)$. In fact, the restriction of $x(\cdot, t)$ on $\p M$ is contained in $\RR^n$ and hence,
$$f\nu+T|_{\p \S_t}= \p_t x|_{\p M}\in T\RR^n.$$
From \eqref{co-normal bundle}, we know 
$$\nu=\frac{1}{\sin \theta} \bar \nu -\cot\theta \mu.$$ 
Since $\bar \nu\in T\RR^n$,  it follows $T-f\cot\theta \mu \in  T\RR^n\cap T\S_t=T(\p \Sigma_t)$,
 and hence \eqref{choiceT1}. 
Up to a diffeomorphism of $\p  M$, we can assume $\widetilde{T}=0$. 
For simplicity, in the following, we always assume that 
\begin{eqnarray}\label{choiceT}
	T |_{\p \S_t}=f\cot\theta \mu.
\end{eqnarray}
Hence, from now on, let $\Sigma_t$ be a family of smooth, embedding hypersurfaces with $\theta$-capillary boundary in $\overline{ \mathbb{R}}^{n+1}_+$, given by the embeddings $x(\cdot,t):M\to \overline{\mathbb{R}}^{n+1}_+$, which evolves by the general flow
\begin{eqnarray}\label{flow with normal and tangential}
	\p_t x=f\nu+T,
\end{eqnarray} with  $T\in T\Sigma_t$  satisfying    \eqref{choiceT}.   {We emphasize that the tangential part $T$ plays a key role in the proof of Theorem \ref{1st variation of total sigma k} below.
}

Along flow \eqref{flow with normal and tangential}, we have the following evolution equations for the induced metric $g_{ij}$, the area element $dA_t$, the
unit outward normal $\nu$, the second fundamental form $h_{ij}$, the Weingarten
matrix  $h^i_j$,  the mean curvature $H$, the $k$-th mean curvature $\sigma_k$ and $F:=F(h^j_i)$ of the hypersurfaces $\Sigma_t$.  These evolution equations will be used later.

\begin{prop}\label{basic evolution eqs}
	Along  flow \eqref{flow with normal and tangential}, it holds that
	\begin{enumerate} 
		\item $\p_t g_{ij}=2fh_{ij}+\n_i T_j+\n_jT_i$.
		\item $\p_tdA_t =\left(fH+\div(T)\right)dA_t.$ 
		\item $\p_t\nu =-\n f+h(e_i,T)e_i$.	
		\item $\p_t h_{ij}=-\n^2_{ij}f +fh_{ik}h_{j}^k +\n_T h_{ij}+h_{j}^k\n_iT_k+h_{i}^k\n_j T_k.$
		\item $\p_t h^i_j=-\n^i\n_{j}f -fh_{j}^kh^{i}_k+\n_T h^i_j.$
		\item $\p_t H=-\Delta f-|h|^2 f+ \langle\n  H, T\rangle $.
		\item $\p_t \s_k=- \frac{\p \s_k}{\p h_i^j}\n^i\n_{j}f - f (\s_1\s_k-(k+1)\s_{k+1})+\langle \n \s_k, T\rangle $.
		\item  $\p_t F= -F^{j}_i\n^i\n_{j}f -fF^{j}_ih_{j}^kh^{i}_k+\langle\n F,T\rangle $,  where $F^i_j:=\frac{\p F}{\p h^j_i}$.
		 
	\end{enumerate}
\end{prop}  
The proof of Proposition \ref{basic evolution eqs} for $T=0$ can be found for example in  \cite[Chapter 2, Section 2.3]{Ger} or \cite[Appendix B]{Ec}. 
	A proof for a general $T$ can be found in \cite[Proposition 2.11]{WeX21}.

Now  we complete the proof of Theorem \ref{1st variation of total sigma k}.
\begin{proof}[\textbf{Proof of Theorem \ref{1st variation of total sigma k}}]
Choose an orthonormal frame $\{e_\alpha\}_{\alpha=2}^{n}$ of $T{\p\S}\subset T\mathbb{R}^n$ such that $\{e_1:=\mu,(e_\alpha)_{\alpha=2}^n\}$ forms an orthonormal frames for $T\S$. 
	First, by taking time derivative to the capillary boundary condition, $\langle \nu,\overline{N}\circ x\rangle=-\cos\theta$ along $\p\S$, we obtain
	\begin{eqnarray*}
		0&=&\langle \p_t \nu, \overline{N}(x(\cdot,t)) \rangle +\langle \nu, d\overline{N}(f\nu+T)\rangle
		\\&=&\langle -\n  f+h(e_i,T)e_i,\overline{N}\rangle  
 	\\&=&-\sin\theta \n_{\mu} f+\sin\theta h(e_i,T) \langle e_i,\mu\rangle 
		\\&=&-\sin\theta \n_\mu f+\sin\theta h(\mu,\mu) \cot \theta  f ,
	\end{eqnarray*}where we have used  \eqref{co-normal bundle},  Proposition \ref{basic evolution eqs} and $T|_{\p M}=f\cot\theta\mu$. As a result, 
	\begin{eqnarray}\label{robin bdry of speed}
		\n_\mu f= \cot\theta h(\mu,\mu) f \text{ on }\p \S_t.
	\end{eqnarray}
Next, using integration by parts and Proposition \ref{basic evolution eqs} we have
\begin{eqnarray}
\frac{d}{dt} \left(\int_{\Sigma_t} \sigma_k dA_t\right) &=&\int_{\Sigma_t} \big[(\p_t \sigma_k ) dA_t+\sigma_k \p_t (dA_t)\big] \notag
	\\&=&\int_{\Sigma_t} \frac{\p\sigma_k}{\p h_i^j}\big( -f^i_{;j}-fh_{ik}h^{kj}+\langle \n  h^i_j, T\rangle\big)   dA_t+\int_\S \sigma_k\big(f\sigma_1+\div(T)\big)dA_t\notag
	\\&=&-\int_{\p\S_t} \frac{\p\sigma_{k}}{\p h^j_i} f^i\mu_j+\int_{\p\S_t} \sigma_k \langle T,\mu\rangle +\int_{\Sigma_t} f\big( \sigma_1\sigma_k-\frac{\p\sigma_k}{\p h^j_i} h_{ik}h^{kj}\big)dA_t\notag
	\\&=&\int_{\p\S_t} \big(f\sigma_k \cot\theta-\frac{\p\sigma_{k}}{\p h^j_i} f^i\mu_j\big) +(k+1)\int_{\Sigma_t} f\sigma_{k+1}dA_t\notag
	\\&=&\int_{\p\Sigma_t} \cot\theta f   \sigma_k (h|h_{11})    +(k+1)\int_{\Sigma_t} f\sigma_{k+1}dA_t,\label{first var of total sigma k}
\end{eqnarray}where we have used $T|_{\p M}=f\cot \theta \mu$,    \eqref{robin bdry of speed} and Lemma \ref{prop2.2} (1), (4).

Moreover  flow \eqref{flow with capillary} induces a hypersurface flow $\p\S_t\subset \RR^n$ with normal speed $\frac{f}{\sin\theta}$, that is,
\begin{eqnarray*}
	\p_t x|_{\p M}=f\nu+f\cot\theta \mu=\frac{f}{\sin\theta}\overline{\nu}.
\end{eqnarray*}By  \eqref{add_1}, we have
\begin{eqnarray*}
	\frac{d}{dt}\V_k(\widehat{\p\S_t})=\frac{n-k}{n} \int_{\p\S_t} \frac{f}{\sin\theta} H_k^{\p\S_t}(\widehat{h}).
\end{eqnarray*}From Proposition \ref{basic-capillary} (2), we know
\begin{eqnarray*}
h_{\alpha\beta} = \sin\theta \widehat{h}_{\alpha\beta},
\end{eqnarray*}
and hence $\sigma_k(h|h_{11})=\sin^k\theta \sigma_k(\widehat{h})$.  Substituting these formulas into \eqref{first var of total sigma k}, we obtain
\begin{eqnarray*}
\frac{d}{dt} \big(\int_{\S_t} H_k dA_t- \sin^k\theta \cos\theta \V_k(\widehat{\p\S_t})\big)=(n-k)\int_{\S_t} f H_{k+1} dA_t.
\end{eqnarray*}By the definition of $\V_{k+1,\theta}(\widehat{\S_t})$ in \eqref{quermassintegrals}, we get the desired formula
\eqref{fvf k+1} for $k\ge 0.$

It remains to consider the case  $k=-1$. It is easy to check that
$$\V_{0,\theta}(\Sigma_t) =|\widehat \Sigma_t| =\frac 1{n+1}\int_{\Sigma_t} \< x, \nu \> dA_t.$$ A direct computation gives
\begin{eqnarray*}
(n+1) \frac d {dt} \V_{0,\theta}(\Sigma_t) &=&  \int_{\S_t}  \left[f - \<x, \n f  \> + \<x,\nu\>  fH +
h(T, x^T) + \<x, \nu\>\div \, T \right]dA_t
\\  &=& 
 \int_{\S_t} \left(  (1 + \div (x^T) )f   + \<x,\nu\>  fH \right) dA_t + \int_{\p \S_t} \left(-\<x^T, \mu\>f+ \<x,\nu \> \< T,\mu\> \right) \\  &=& 
(n+1) \int _{\S_t} f dA_t, 
\end{eqnarray*} since $-\<x^T, \mu\> f +\<x,\nu \> \< T, \mu\> =0$ for $x\in \p\S$, which
 follows from 
$$\<x,\nu \> \< T, \mu\> = f\cot \theta \<x,\nu \>= f\cos \theta \<x,\bar \nu\>,
$$
$$\<x^T, \mu\> f = f\cos \theta \<x,\bar \nu\>.
$$
Now we complete the proof. \end{proof}

\section{Locally constrained curvature flow}\label{sect 3}
In this section,  we first introduce a new locally constrained curvature flow and show the monotonicity of the quermassintegral along the flow.

Let $M$ be a compact orientable smooth $n$-dimensional manifold. Suppose $x_0:M\to \overline{\R}^{n+1}_+$ be a smooth initial embedding such that $x_0(M)$ is a convex hypersurface  in $\overline{\RR}^{n+1}_+$ and intersects with $\p{\RR}^{n+1}_+$ at a constant contact angle $\theta\in (0,{\pi})$. 
We consider the smooth family of   embeddings    $x:M\times [0,T)\to \overline{\mathbb{R}}^{n+1}_+, $ satisfying the following evolution equations
\begin{equation}\label{flow with capillary}
	\begin{array}{rcll}
	( \partial_t x(p,t)  )^\perp &=&f(p,t)\nu \quad &
		\hbox{ for } (p,t)\in M\times[0,T),\\
		\langle \nu(p,t),\overline{N}\circ x(p,t)\rangle  &=&  \cos(\pi-\theta )
		\quad & \hbox{ for  }(p,t)\in\partial M\times [0,T),
	\end{array}
\end{equation}with $x(M,0)=x_0(M)$ and
 \begin{eqnarray}\label{speed}
	f:=\frac{1+\cos\theta \langle \nu,e\rangle }{F}-\langle x,\nu\rangle ,
\end{eqnarray} 
where
\begin{eqnarray}\label{f0}
F:= \frac{H_{l} }{H_{l-1} }.
\end{eqnarray}

\
\
The following nice property of  flow \eqref{flow with capillary}  is essential for us to prove Theorem \ref{thm 1} later.
\begin{prop}\label{monotone along flow}
	As long as  flow \eqref{flow with capillary} exists and $ \S_t$ is strictly $l$-convex, 
	$ \V_{l,\theta}(\widehat{\S_t})$ is preserved and $ \V_{k,\theta}(\widehat{\S_t})$ is non-decreasing  for $1\leq k<l \leq n$.
\end{prop}
\begin{proof}
From Theorem \ref{1st variation of total sigma k}, we see
	\begin{eqnarray*}
		\p_t \V_{l,\theta}(\widehat{\S_t})&=&\frac{n+1-l}{n+1}\int_{\S_t} f  H_l dA_t\\&=&\frac{n+1-l}{n+1}\int_{\S_t}\left[ \left(1+\cos\theta \langle \nu,e\rangle \right)H_{l-1}-H_l \langle x,\nu\rangle\right]dA_t\\& =&0,
	\end{eqnarray*} where the last equality follows from \eqref{minkowski sigma_k}. 
For $1\leq k<l\leq  n$, from Theorem \ref{1st variation of total sigma k}
	\begin{eqnarray*}
		\p_t \V_{k ,\theta}(\widehat{\S_t})&=&\frac{n+1-k}{n+1}\int_{\S_t} f H_{k } dA_t
		\\&=& \frac{n+1-k}{n+1}\int_{\S_t} \left[H_{k}\frac {H_{l-1}}{{H_l}}\big(1+\cos\theta \langle \nu,e\rangle\big)-H_k\langle x,\nu\rangle\right]dA_t
		\\&\geq& \frac{n+1-k}{n+1}\int_{\S_t} \left[H_{k-1}\big(1+\cos\theta \langle \nu,e\rangle\big)-H_k\langle x,\nu\rangle\right]dA_t
		\\&=&0,
	\end{eqnarray*} where we have used the Newton-MacLaurin inequality \eqref{NM-ineq} and the  Minkowski formula \eqref{minkowski sigma_k} in the last two steps respectively.
 \end{proof}


\

\section{A priori estimates and convergence}\label{sect 4}
The main result of this section is the following  long-time existence and  the convergence result of  flow \eqref{flow with capillary}
with $l=n$, i.e.,
\begin{eqnarray}\label{eq4.1}
F=\frac{H_n}{H_{n-1}}\end{eqnarray} 
under an angle constraint
$$\theta \in (0,\frac{\pi}{2}].$$

  \begin{thm}\label{longtime exist and conv} Assume $x_0:M\to \overline{\mathbb{R}}^{n+1}_+$ is
 an embedding of a strictly convex  capillary hypersurface 
  in the half-space with the contact angle $\theta \in ( 0,\frac{\pi}{2}]$. 
  Then there exists $x:M\times[0,+\infty)\to \overline{\mathbb{R}}^{n+1}_+$ satisfying flow \eqref{flow with capillary}
  with $F$ given by \eqref{eq4.1} and the initial condition  $x(M,0)=x_0(M)$. Moreover, $x(\cdot,t)\to x_\infty(\cdot)$ in $C^\infty$ topology as $t\to+\infty$, and the limit $x_\infty:M\to \overline{\mathbb{R}}^{n+1}_+$ is a spherical cap. 
 \end{thm}
In order to prove this theorem, we need to obtain  a priori estimates, which will be given as follows.
\subsection{The short time existence}
For the short time existence, one can follow the strategy presented in the paper of Huisken-Polden \cite{HP} to give a proof for a general initial capillary hypersurface. Since our initial hypersurface is convex, one can prove the short time existence
in the class of star-shaped hypersurfaces. In this class, one can in fact reduce  flow \eqref{flow with capillary} to a scalar flow. Then the short time existence follows clearly from the standard theory for parabolic equations.
Therefore we first consider the reduction.

Assume that  a capillary hypersurfaces 
 $\S$ is  strictly star-shaped with respect to the origin. 
 One can reparametrize it  as a graph over $\overline{\mathbb{S}}^n_+$. Namely, 
 there exists a positive function $r$ defined on $\overline{\SS}^n_+$ such that
 \begin{eqnarray*}
 	\S=\left\{ r(X )X |    X\in \overline{\SS}^n_+\right\},
 \end{eqnarray*}where $X:=(X_1,\ldots, X_n)$ is a local coordinate of $\overline{\SS}^n_+$. 
 
 We  denote  $\n^0$ be the Levi-Civita connection on $\SS^n_+$ with respect to the standard round metric $\sigma:=g_{_{\SS^n_+}}$, $\p_i:=\p_{X_i}$,  $\sigma_{ij}:=\sigma(\p_{i},\p_{j})$,  $r_i:=\n^0_{i} r$, and $r_{ij}:=\n^0_{i}\n^0_j r$. The induced metric $g$ on $\S$ is given by
 \begin{eqnarray*}
 	g_{ij}=r^2\sigma_{ij}+ r_ir_j=e^{2\varphi}\left(\sigma_{ij}+\varphi_i\varphi_j\right),
 \end{eqnarray*}where  $\varphi(X):=\log r(X)$. Its inverse $g^{-1}$ is given by
 \begin{eqnarray*}
 	g^{ij}=\frac{1}{r^2} \left(\sigma^{ij}-\frac{r^ir^j}{r^2+|\n^0 r|^2}\right)=e^{-2\varphi}\left( \sigma^{ij}-\frac{\varphi^i\varphi^j}{v^2}\right),
 \end{eqnarray*}where $r^i:=\sigma^{ij}r_j$,  $\varphi^i:=\sigma^{ij}\varphi_j$ and $v:=\sqrt{1+|\n^0 \varphi|^2}$.
 The unit outward normal vector field  on $\S$ is given by
 \begin{eqnarray*}
 	\nu=\frac{1}{v}\left( \p_r-r^{-2}\n^0 r\right)=\frac{1}{v}\left( \p_r-r^{-1}\n^0 \varphi\right).
 \end{eqnarray*}
 The second fundamental form $h$ on $\S$ is
 \begin{eqnarray*}
 	h_{ij}=\frac{e^\varphi}{v }\left( \sigma_{ij}+\varphi_i\varphi_j-\varphi_{ij}\right),
 \end{eqnarray*}and its Weingarten matrix $h^i_j$ is
 \begin{eqnarray*}
 	h_j^{i}=g^{ik}h_{kj}=\frac{1}{e^\varphi v }\left[ \delta^i_j-(\sigma^{ik}-\frac{\varphi^i\varphi^k}{v^2})\varphi_{kj}\right].
 \end{eqnarray*}
The higher order mean curvature $H_k$ can also be expressed by $\varphi$.
Moreover, 
 \begin{eqnarray*}
 	\langle x,\nu\rangle =\langle r\p_r,\nu\rangle =\frac{e^\varphi}{v }.
 \end{eqnarray*}
 
 In order to express the capillary boundary condition in terms of the radial function $\varphi$, we use the polar coordinate in 
 the half-space. For $x:=(x',x_{n+1})\in \RR^n\times [0,+\infty)$  and $X:=(\beta,\xi)\in [0,\frac{\pi}{2}]\times \SS^{n-1}$, we have that 
 \begin{eqnarray*}
 	x_{n+1}=r\cos\beta,\quad |x'|=r\sin\beta.
 \end{eqnarray*} 
  Then
 \begin{eqnarray*}
 	e_{n+1}=	\p_{x_{n+1}}=\cos\beta \p_r-\frac{\sin\beta}{r}\p_\beta.
 \end{eqnarray*} 
 In these coordinates the standard Euclidean metric is given by
 $$ |dx|^2=dr^2+r^2\left(d\beta^2+\sin^2\beta g_{_{\SS^{n-1}}}\right).$$
 It follows that
 \begin{eqnarray*}
 	\langle \nu,e_{n+1}\rangle =\frac{1}{v}\left(\cos\beta+\sin\beta \n_{\p_\beta}^0   \varphi \right).
 \end{eqnarray*}

 Along $\p\SS^n_+$ it holds
 \begin{eqnarray*}  	\overline{N}\circ x=-e_{n+1}
	=\frac{1}{r}\p_\beta,
 \end{eqnarray*}
 which yields
 \begin{eqnarray*}
 	-\cos\theta =	\langle \nu,\overline{N}\circ x\rangle =\langle \frac{1}{v}\left( \p_r-r^{-1}\n^0 \varphi\right), \frac{1}{r}\p_\beta\rangle=-\frac{\n^0_{\p_\beta}\varphi}{v} ,
 \end{eqnarray*} that is,
 \begin{eqnarray} \label{4.1b}
 	\n^0_{\p_\beta} \varphi=\cos\theta  \sqrt{1+|\n^0\varphi|^2}.
 \end{eqnarray}
 Therefore, in the class of star-shaped hypersurfaces flow  \eqref{flow with capillary} is reduced to
  the following scalar  parabolic equation with an oblique boundary condition
 \begin{eqnarray}\label{scalar flow with capillary}
 	\begin{array}{rcll}
 		\p_t \varphi &=& \frac{v}{e^\varphi}f,\quad &\text{ in } \SS^n_+\times [0,T^*),  \\
 		\n_{\p_\beta}^0 \varphi&=&\cos\theta \sqrt{1+|\n^0\varphi|^2},\quad &\text{ on }
		\p\SS^n_+\times [0,T^*),\\
 		\varphi(\cdot,0)&=&\varphi_0(\cdot),\quad &\text{ on } \SS^n_+,
 	\end{array}
 \end{eqnarray}
 where $\varphi_0 $ is the parameterization radial function of $x_0(M)$ over $\overline{\SS}^n_+$, and
 \begin{eqnarray*}
 	f=\frac{H_{n-1}}{H_n}\left[1- \frac{\cos\theta}{v}\left( \cos\beta+\sin\beta \n^0_{\p_\beta}\varphi \right)  \right]-\frac{e^\varphi}{v}.
 \end{eqnarray*}
Since $|\cos \theta|<1$, the oblique boundary condition \eqref{4.1b} satisfies the non-degeneracy condition
in \cite{NU}, see also \cite{Dong}.
 Hence the short time existence follows.

\subsection{Barriers}
Let $T^*$ be the maximal time of smooth existence of a solution to \eqref{flow with capillary}, more precisely
in the class of star-shaped hypersurfaces. It is obvious that $F$ can not be zero and hence $F$ is positive in 
$M\times [0, T^*)$.
 The positivity of $F$ implies that $\S_t$  is strictly convex up to $T^*$.

The convexity of $\S_0$ implies that there exists some $0<r_1<r_2<\infty$, such that 
$$\S_0\subset \widehat{C_{r_2,\theta}}\setminus \widehat{C_{r_1,\theta}}.$$
The family of $C_{ r,\theta}$ forms natural barriers of \eqref{flow with capillary}. Therefore, we can show that the solution to \eqref{scalar flow with capillary} is uniformly bounded from above and below.
\begin{prop}\label{barrier}For any $ t\in [0, T^*)$, $\S_t$ satisfies
	$$\S_t\subset \widehat{C_{r_2,\theta} }\setminus \widehat{C_{r_1,\theta} }.$$
\end{prop}
 \begin{proof}
 	Recall that $C_{r,\theta}$ satisfies \eqref{static model radius}. Thus for each $r>0$, it is a static solution to  flow \eqref{flow with capillary}. The assertion follows from the avoidance principle for strictly parabolic equation with a capillary boundary condition (see \cite[Section 2.6]{AW} or \cite[Proposition 4.2]{WW}).  
 \end{proof}

 \subsection{Evolution equations of $F$ and $H$}
We first introduce a   parabolic operator for \eqref{flow with capillary}
\begin{eqnarray*}
	\L :=\p_t-\frac{1+\cos\theta\langle \nu,e\rangle}{F^2}F^{ij}\n_{ij}^2-\langle T+x-\frac{\cos\theta}{F}e,\n\rangle.
\end{eqnarray*}
 Set $\mathcal{F}:=\sum_{i=1}^n F^{i}_i$.  Using Proposition \ref{prop2.2}
 we have
 \begin{eqnarray}
 \mathcal{F}-\frac{F^{ij}h_{ij}}{F}&=& \mathcal{F}-1\geq 0,\label{xeq11}\\
 \frac{F^{ij}h_i^k h_{kj}}{F^2}&=& 1.\label{xeq22}
 \end{eqnarray}

\begin{prop}\label{ev eq of F} Along  flow \eqref{flow with capillary}, we have
	\begin{eqnarray*}
		\L F=2\cos\theta F^{-2}F^{ij}F_{;j} h_{ik}\langle e_k,e\rangle -2(1+\cos\theta \langle \nu,e\rangle)F^{-3}F^{ij}F_{;i}F_{;j}+F\left(1-\frac{F^{ij}(h^2)_{ij}}{F^2}\right),
	\end{eqnarray*}
and  
\begin{eqnarray}\label{neumann of F}
	\n_\mu F=0,\quad \text{ on } \p\S_t.
\end{eqnarray}
\end{prop}
 
\begin{proof}Using the Codazzi formula, we have
	\begin{eqnarray*}
		F^{ij}\langle x,\nu\rangle_{;ij}&=F^{ij}\left(h_{ij}+h_{ij;k}\langle x,e_k\rangle -(h^2)_{ij}\langle x,\nu\rangle \right)\\&=F+F_{;k}\langle x,e_k\rangle -F^{ij}(h^2)_{ij}\langle x,\nu\rangle ,
	\end{eqnarray*}
and
\begin{eqnarray*}
	F^{ij}\langle \nu,e\rangle_{;ij}&=&F^{ij}\big(h_{ik;j}\langle e_k,e\rangle -(h^2)_{ij}\langle \nu,e\rangle\big)
	\\&=& F_{;k}\langle e_k,e\rangle -F^{ij}(h^2)_{ij}\langle \nu,e\rangle .
\end{eqnarray*}
Combining with Proposition \ref{basic evolution eqs}, we obtain
	\begin{eqnarray*}
		\p_t F&=&-F^{ij}f_{;ij}-fF^{ij}(h^2)_{ij}+\langle \n F,T\rangle
		\\&=&-F^{ij} \big( \frac{1+\cos\theta \langle \nu,e\rangle }{F}-\langle x,\nu\rangle  \big)_{;ij}-fF^{ij}(h^2)_{ij}+\langle \n F,T\rangle
		\\&=&-\cos\theta F^{ij}F^{-1} \langle \nu,e\rangle_{;ij}+2\cos\theta F^{-2}F^{ij}F_{;j} \langle \nu,e\rangle_{;i}-2(1+\cos\theta \langle \nu,e\rangle)F^{-3}F^{ij}F_{;i}F_{;j}\\&&+F^{-2}F^{ij}F_{;ij}(1+\cos\theta \langle \nu,e\rangle)+\big(F+F_{;k}\langle x,e_k\rangle -F^{ij}(h^2)_{ij}\langle x,\nu\rangle \big)
		\\&& -(1+\cos\theta\langle \nu,e\rangle) F^{-1}F^{ij}(h^2)_{ij}+\langle x,\nu\rangle F^{ij}(h^2)_{ij}+\langle \n F,T\rangle.
	\end{eqnarray*} 
Hence it follows
\begin{eqnarray*}
\mathcal{L}F&=&\p_tF -(1+\cos\theta \langle \nu,e\rangle)F^{-2}F^{ij}F_{;ij}-\langle T+x-\cos\theta F^{-1}e,\n F\rangle
\\&=&  \cos\theta F^{ij}(h^2)_{ij}F^{-1} \langle \nu,e\rangle +2\cos\theta F^{-2}F^{ij}F_{;j} h_{ik}\langle e_k,e\rangle -2(1+\cos\theta \langle \nu,e\rangle)F^{-3}F^{ij}F_{;i}F_{;j}\\&&+\big(F -F^{ij}(h^2)_{ij}\langle x,\nu\rangle \big)
 -(1+\cos\theta\langle \nu,e\rangle) F^{-1}F^{ij}(h^2)_{ij}+\langle x,\nu\rangle F^{ij}(h^2)_{ij}
 \\&=&   2\cos\theta F^{-2}F^{ij}F_{;j} h_{ik}\langle e_k,e\rangle -2(1+\cos\theta \langle \nu,e\rangle)F^{-3}F^{ij}F_{;i}F_{;j}\\&&+ F  
 -  F^{-1}F^{ij}(h^2)_{ij} .
\end{eqnarray*}

Along $\p\S_t,$ from \eqref{robin bdry of speed} we know  
\begin{eqnarray*}
	\n_\mu f=\cot\theta h(\mu,\mu) f, 
\end{eqnarray*}
By  \eqref{co-normal bundle} or \eqref{2.5}  and {Proposition \ref{basic-capillary} (1)},  
we have on $\p \Sigma_t$ \begin{eqnarray*} 
	\n_\mu \langle x,\nu\rangle  =\langle x, h(\mu,\mu)\mu\rangle =\cos \theta h(\mu,\mu)\langle x, \overline \nu\rangle =  \cot\theta h(\mu,\mu)\langle x,  \nu\rangle 
	 ,
\end{eqnarray*} and hence 
\begin{eqnarray*}
	\n_\mu (f+\langle x,\nu\rangle )=\cot\theta h(\mu,\mu) (f+\langle x,\nu\rangle ).
\end{eqnarray*}
Using  \eqref{co-normal bundle} and {Proposition \ref{basic-capillary} (1) again}, we have ($e=\overline N$)
\begin{eqnarray*} 
		\n_\mu \langle \nu,e\rangle =h(\mu,\mu)\langle \mu,e\rangle 
	  =-\tan \theta h(\mu,\mu)\langle \nu,e\rangle
\end{eqnarray*}
and
\begin{eqnarray*}
	\n_\mu\big(1+\cos\theta \langle \nu,e \rangle \big)= -\sin\theta h(\mu,\mu)\langle \nu,e\rangle,
\end{eqnarray*}where we used $\langle e,\overline{\nu}\rangle =0$ and $\langle e,x\rangle =0$ on $\p\S_t$.
One can easily check that the left hand side of the previous formula equals to 
$\cot \theta h(\mu,\mu) (1+\cos\theta \langle \nu, e \rangle )$, on $\p \Sigma$.
Hence it follows that
\begin{eqnarray*}
	\n_\mu F&=&\n_\mu \big(\frac{1+\cos\theta \langle \nu,e\rangle }{f+\langle x,\nu\rangle }\big) =0. 
\end{eqnarray*}

\end{proof}
We remark that \eqref{neumann of F} plays an important role in applying the maximum principle later. This property  holds for curvature flow of free boundary hypersurfaces
and capillary hypersurfaces,  see also \cite{SWX, WeX21}.
\begin{prop}\label{ev eq of H} Along flow \eqref{flow with capillary}, we have
	\begin{eqnarray} \label{evo_H}
		\L H&=&(1+\cos\theta \langle \nu,e\rangle )F^{-2} F^{kl,st}h_{kl;i}h_{st;i}   +(2+\cos\theta \langle \nu,e\rangle)   H  \nonumber \\&& 
		+\Big[ 2\cos\theta F^{-2}F_{;i} \langle \nu,e\rangle_{;i} -2\left(1+\cos\theta \langle \nu,e\rangle\right) F^{-3}|\n F|^2\\&&-(2+\cos\theta \langle \nu,e\rangle) F^{-1}|h|^2  \Big],
	\nonumber \end{eqnarray}
	and, {while $\S$ is convex,}
	\begin{eqnarray}\label{neumann of H}
		\n_\mu H\leq 0,\quad \text{ on } \p\S_t.
	\end{eqnarray}
\end{prop}

\begin{proof}
First, note that
	\begin{eqnarray*}
		\Delta \langle \nu,e\rangle =H_{;k}\langle e_k,e\rangle -|h|^2\langle \nu,e\rangle,
	\end{eqnarray*}
	\begin{eqnarray*}
	\Delta \langle x,\nu\rangle =H+H_{;k}\langle x,e_k\rangle -|h|^2\langle x,\nu\rangle .
\end{eqnarray*}
Applying  Proposition \ref{basic evolution eqs}, we obtain 
	\begin{eqnarray*}
		\p_t H&=&-\Delta f-|h|^2 f+\langle  \n H,T\rangle 
		\\&=& (1+\cos\theta \langle \nu,e\rangle)F^{-2}\Delta F-2(1+\cos\theta \langle \nu,e\rangle ) F^{-3}|\n F|^2 \\&& +2\cos\theta F^{-2}F_{;i} \langle \nu,e\rangle_{;i}-\cos\theta F^{-1}\Delta \langle \nu,e\rangle +\Delta \langle x,\nu\rangle 
		\\&&-(1+\cos\theta \langle \nu,e\rangle )F^{-1}|h|^2+\langle x,\nu\rangle |h|^2+\langle  \n H,T\rangle 
		\\&=&(1+\cos\theta \langle \nu,e\rangle)F^{-2}\Delta F-2(1+\cos\theta \langle \nu,e\rangle ) F^{-3}|\n F|^2 \\&& +2\cos\theta F^{-2}F_{;i} \langle \nu,e\rangle_{;i}- F^{-1}|h|^2	
	+H\\&&+ \langle x,\n H\rangle 
 +\langle  \n H,T\rangle  -\cos\theta F^{-1}\langle e,\n H\rangle.
	\end{eqnarray*}
The Ricci equation and the Codazzi   equation  yield
\begin{eqnarray*}
	h_{kl;ii}&=&h_{ki;li}=h_{ki;il}+R^p_{ili}h_{pk}+R^p_{kli}h_{pi}
	\\&=&h_{ii;kl}+(h_{pl}H-h_{pi}h_{li})h_{pk}+(h_{pl}h_{ki}-h_{pi}h_{kl})h_{pi}
	\\&=&H_{;kl}+h_{pk}h_{pl}H-|h|^2 h_{kl},
\end{eqnarray*}which implies 
 \begin{eqnarray*}
	\Delta F&=&\frac{\p^2 F}{\p h_{kl}\p h_{st}} h_{kl;i}h_{st;i}+F^{kl}h_{kl;ii}
 \\&=&F^{kl,st}h_{kl;i}h_{st;i}+F^{kl} H_{;kl}+F^{kl}(h^2)_{kl}H-F|h|^2  .
\end{eqnarray*}
Hence we have
	\begin{eqnarray*}
\L H&=&\p_t H-F^{-2}(1+\cos\theta \langle \nu,e\rangle )F^{ij}H_{;ij}-\langle T+x-\cos\theta F^{-1}e,\n H\rangle 
	\\&=&(1+\cos\theta \langle \nu,e\rangle)F^{-2}\big[ F^{kl,st}h_{kl;i}h_{st;i} +F^{kl}(h^2)_{kl}H-F|h|^2\big]\\&&-2(1+\cos\theta \langle \nu,e\rangle ) F^{-3}|\n F|^2  +2\cos\theta F^{-2}F_{;i} \langle \nu,e\rangle_{;i}	- F^{-1}|h|^2	
+H
\\&=&(1+\cos\theta \langle \nu,e\rangle )F^{-2} F^{kl,st}h_{kl;i}h_{st;i}   +(2+\cos\theta \langle \nu,e\rangle)   H\\&& 
+\Big[ 2\cos\theta F^{-2}F_{;i} \langle \nu,e\rangle_{;i} -2(1+\cos\theta \langle \nu,a\rangle) F^{-3}|\n F|^2\\&&-(2+\cos\theta \langle \nu,e\rangle) F^{-1}|h|^2  \Big].
\end{eqnarray*}

Along $\p\S_t$, choosing an orthonormal frame $\{e_\alpha\}_{\alpha=2}^{n}$ of $T{\p\S_t}$ such that $\{e_1:=\mu,(e_\alpha)_{\alpha=2}^n\}$ forms an orthonormal frames for $T\S_t$. From Proposition \ref{basic-capillary}, we have
\begin{eqnarray*}
	h_{\alpha\beta;\mu}= \cos\theta \widehat{h}_{\beta\gamma} (h_{11}\delta_{\alpha\gamma}-h_{\alpha\gamma}),
\end{eqnarray*}for all $2\leq \alpha\leq n$. 
Equation \eqref{neumann of F} implies $$0=\n_\mu F=F^{11}h_{11;1}+\sum\limits_{\alpha=2}^nF^{\alpha\alpha}h_{\alpha\alpha;1},$$ 
which in turn implies  
\begin{eqnarray*}
	\n_\mu H&=&h_{11;1}+\sum_{\alpha=2}^n h_{\alpha\alpha;1}
	\\&=&-\sum_{\alpha=2}^n \frac{F^{\alpha\alpha}}{F^{11}}h_{\alpha\alpha;1}+\sum_{\alpha=2}^n h_{\alpha\alpha;1}
	=\sum_{\alpha=2}^n \frac{1}{F^{11}}\big( F^{11}-F^{\alpha\alpha}\big)h_{\alpha\alpha;1}
	\\&=&\sum_{\alpha=2}^n \frac{1}{F^{11}}\big( F^{11}-F^{\alpha\alpha}\big)(h_{11}-h_{\alpha\alpha})  \widetilde{h}_{\alpha\alpha} 
	\\&\leq& 0,
\end{eqnarray*}where the last inequality follows from the concavity of $F$, and the  convexity of $\p\S\subset \S$, {see Corollary \ref{x-conv}}. Hence \eqref{neumann of H} is proved.
	\end{proof}
{\begin{rem}
\eqref{neumann of H} is the only place where we have used $\theta\in (0,\frac{\pi}{2}]$.
\end{rem}	
}
\subsection{Curvature estimates}\
 First, we have the uniform  bound of $F$, which follows directly from Proposition \ref{ev eq of F} and the maximum principle.
\begin{prop}\label{bounds of F}
	Along  flow \eqref{flow with capillary}, it holds
	\begin{eqnarray*}
\min_{M} F(\cdot, 0)\leq		F(p ,t)\leq \max_{M} F(\cdot, 0), \quad \forall (p,t)\in  M\times [0,T^*).
	\end{eqnarray*}
\end{prop}

 In particular, from the uniform lower bound of $F:=\frac{H_n}{H_{n-1}}$, we get a uniform curvature positive lower bound.  \begin{cor}\label{convexity preserve}
 	$\Sigma_t, t\in [0,T^*)$ is uniformly convex, that is, there exists $c>0$ depending only on $\S_0$, such that the principal curvatures of $\S_t$, $$\min_i \kappa_i(p, t) \ge c,$$  for all $(p,t)\in M\times [0,T^*)$.
 \end{cor}

Next we obtain the uniform bound of the mean curvature.
\begin{prop}\label{mean-curv-bound}
There exists $C>0$ depending only on $\S_0$, such that
	\begin{eqnarray*} 
		H(p,t)\leq C,\quad \forall (p,t)\in M\times [0,T^*).
	\end{eqnarray*} 
\end{prop}	
\begin{proof}
From \eqref{neumann of H}, we know that $\n_\mu H\leq 0$ on $\p\S_t$. Thus $H$ attains its  maximum value at some interior point, say $p_0\in \text{int}(M)$. We now compute at $p_0$. 
	
	From the concavity of $F=\frac{n\s_n}{\s_{n-1}}$ in Proposition  \ref{prop2.3}, we know
	$$   \left(1+\cos\theta \langle\nu ,e\rangle\right) F^{-2} F^{kl,st}h_{kl;i}h_{st;i} \le 0.$$
Using Proposition \ref{ev eq of H}, we have
		\begin{eqnarray*}
		\L H&\leq &  (2+\cos\theta \langle \nu,e\rangle)   H 
		\\&&+\Big[ 2\cos\theta F^{-2}F_{;i} \langle \nu,e\rangle_{;i} -2\left(1+\cos\theta \langle \nu,e\rangle\right) F^{-3}|\n F|^2-(2+\cos\theta \langle \nu,e\rangle) F^{-1}|h|^2  \Big]
		\\&:=& \text{K}_1+\text{K}_2,
	\end{eqnarray*}
The term $|\text{K}_1|$ is bounded by $3H$. 
For the  term $\text{K}_2$, we note that
	\begin{eqnarray*}
F\text{K}_2&:=&2\cos\theta F^{-1}F_{;i} \langle \nu,e\rangle_{;i} -2(1+\cos\theta \langle \nu,e\rangle) F^{-2}|\n F|^2 -(2+\cos\theta \langle \nu,e\rangle)  |h|^2 \\&=&2\cos\theta  F^{-1}F_{;i}h_{ii}\langle e,e_i\rangle-2(1+\cos\theta \langle \nu,e\rangle )F^{-2}|\n F|^2-(2+\cos\theta \langle \nu,e\rangle)  |h|^2 
		\\&:=&-\sum_{i=1}^n \big(	\text{S}_1 F_{;i}^2+	\text{S}_{2,i}F_{;i} h_{ii}+	\text{S}_3h_{ii}^2\big)
		=-	\text{S}_1\sum_{i=1}^n{\left(F_{;i}-\frac{	\text{S}_{2,i}}{2	\text{S}_1}h_{ii}\right)^2}+\sum_{i=1}^n \left(\frac{	\text{S}_{2,i}^2}{4	\text{S}_1}-	\text{S}_3\right) h_{ii}^2,
	\end{eqnarray*}
	where we have used the notations 
	\begin{eqnarray*}
		\text{S}_1:=2(1+\cos\theta \langle \nu,e\rangle ) F^{-2}, \quad 	{\text{S}_{2,i}:=-2\cos\theta F^{-1}\<e,e_i\>}, \quad 	\text{S}_3:=2+\cos\theta \langle \nu,e\rangle .
	\end{eqnarray*}
One can check 
	\begin{eqnarray*}
			\text{S}_{2,i}^2-4	\text{S}_1	\text{S}_3&:=& 4\cos^2\theta F^{-2} \langle e,e_i\rangle^2-8(1+\cos\theta \langle \nu,e\rangle ) F^{-2} (2+\cos\theta \langle \nu,e\rangle)
		\\&\leq&  4F^{-2}\Big[ \cos^2\theta |e^T|^2-2(1+\cos\theta \langle \nu,e\rangle )(2+\cos\theta \langle \nu,e\rangle)\Big]
		\\&=&4F^{-2}\Big[-3(1+\cos\theta \langle \nu,e\rangle)^2-1+\cos^2\theta\Big]
		\\&\leq &-c_0
	\end{eqnarray*}for some positive constant $c_0$. Combining with Proposition \ref{bounds of F}, it implies
\begin{eqnarray*}
\text{K}_2\leq -\frac{c_0}{4F \text{S}_1}|h|^2=-\frac{c_0F}{8(1+\cos\theta \<\nu,e\>)} |h|^2\leq -C |h|^2,
\end{eqnarray*} for some positive constant $C>0$. 
Therefore, 
\begin{eqnarray*}
	&&0 \leq \mathcal{L}H(p_0) \leq \text{K}_1+\text{K}_2 \le 3H-C|h|^2,\end{eqnarray*}
which yields that $H$ is uniformly bounded from above.
\end{proof}
 Proposition \ref{mean-curv-bound} and Corollary \ref{convexity preserve} imply  directly that
 \begin{cor}\label{curvature bound}
 	$\Sigma_t, t\in [0,T^*)$, has  a uniform curvature bound, namely,
 	there exists $C>0$ depending only on $\S_0$, such that the principal curvatures of $\S_t$, $$\max_i \kappa_i(p, t) \le C,$$for all $(p,t)\in M\times [0,T^*)$.
 \end{cor}

\subsection{Convergence of the flow}\label{sect 3.3}\

First we show that the convexity implies that the star-shaped is preserved in the following sense. 
 \begin{prop}\label{star-shaped preserve}
 	There exists $c_0>0$ depending only on $\S_0$, such that 
 	\begin{eqnarray}\label{star-shape est}
 		\langle x,\nu\rangle(p,t) \geq c_0.
 	\end{eqnarray}for all $(p,t)\in M\times [0,T^*)$.
 	
 \end{prop}
 
 \begin{proof}For any $T'<T^*$, assume $\min\limits_{M\times [0,T']} \langle x,\nu\rangle(p,t)=\langle x,\nu\rangle (p_0,t_0)$. Then, either $p_0\in\p M$ or $p_0\in M\setminus\p M$. 
 
 		If $p_0\in M\setminus\p M$, let $\{e_i\}_{i=1}^n$ be the orthonormal frame of $\S_t$, then at $p_0$,
 		\begin{eqnarray*}
 			0=	D_{e_i}\langle x,\nu\rangle =h_{ij}\langle x,e_j\rangle.
 		\end{eqnarray*}
 	Due to the strict convexity $(h_{ij})>0$, we  have $\langle x,e_i\rangle=0$. It follows  
 		\begin{eqnarray*}
 			\langle x,\nu\rangle(p_0)=|x|(p_0)\geq c_0,
 		\end{eqnarray*}for some $c_0>0$, which  depends only on the initial datum.
 		
 		 If $p_0\in \p M$, by \eqref{co-normal bundle}  we have 
 		\begin{eqnarray*}
 			\langle x,\nu\rangle =\langle x,\sin \theta \overline{\nu}-\cos\theta e\rangle =\sin\theta \langle x,\overline{\nu}\rangle .
 		\end{eqnarray*}
 	Hence $\langle x,\overline{\nu}\rangle \big|_{\p M}$ attains its minimum value at $p_0$.
 	As above, choosing $\{e_\alpha\}_{\alpha=2}^n$ be the orthonormal frame of $\p\S_t$ in $\mathbb{R}^n$ such that $e_1=\overline{\nu}$,  we have 
 		\begin{eqnarray*}
 			0=	\n^{\RR^n}_{e_\alpha}\langle x,\overline{\nu}\rangle =\widehat{h}_{\alpha\beta}\langle x,e_\beta\rangle,
 		\end{eqnarray*} 
 		By Proposition \ref{basic-capillary} (2) and Corollary \ref{convexity preserve}, we know $(\widehat{h}_{\alpha\beta})>0$,  and hence  we have $x\parallel \overline{\nu}$ at $p_0$ and 
 		\begin{eqnarray*}
 			\langle x,\overline{\nu}\rangle (p_0) =|x|(p_0)\geq c_0,
 		\end{eqnarray*}for some $c_0>0$, which  depends only  on the initial datum.
 Therefore, we finish the proof of \eqref{star-shape est}.
 \end{proof}
 
\begin{prop}\label{global existence}
	Flow  \eqref{flow with capillary} exists for all time with uniform $C^{\infty}$-estimates.
\end{prop}
\begin{proof}
From Proposition \ref{barrier}, Proposition \ref{star-shaped preserve}, Proposition \ref{convexity preserve} and Corollary \ref{curvature bound}, we see that $\varphi$ is uniformly bounded in $C^2(\mathbb{S}^n_+\times [0,T^*))$ and the scalar equation in \eqref{scalar flow with capillary} is uniformly parabolic.
Since $|\cos\theta|<1$, the boundary value condition in \eqref{scalar flow with capillary} satisfies	the uniformly oblique property. From the standard parabolic theory  (see e.g. \cite[Theorem 6.1, Theorem 6.4 and Theorem 6.5]{Dong}, also  \cite[Theorem 5]{Ura}  and \cite[Theorem 14.23]{Lie}), we conclude the uniform $C^\infty$-estimates and the long-time existence of solution to   \eqref{scalar flow with capillary}.

\end{proof}

\begin{prop}\label{converge limit}
	$x(\cdot, t)$ smoothly converges to a uniquely determined spherical cap around $\cos\theta e$ with capillary boundary, as $t\to\infty$.		
\end{prop}
\begin{proof}
By Proposition \ref{monotone along flow}, we know $ \V_{1,\theta}(\widehat{\S_t})$ is non-decreasing, due to
	\begin{eqnarray*}
		\p_t  \V_{1,\theta}(\widehat{\S_t}) =\frac{n}{n+1}\int_{\S_t} \left(\frac{H_1H_{n-1}}{H_n}
		-1\right) \left(1+\cos\theta  \langle\nu,e\rangle\right)dA_t \geq 0.
	\end{eqnarray*}
{It follows from the  long time existence and  uniform $C^\infty$-estimates that	$$ \int_0^{\infty }	\p_t  \V_{1,\theta}(\widehat{\S_t})  dt \le  \V_{1,\theta}(\widehat{\S_\infty})<+\infty.$$ } Then  we obtain
	\begin{eqnarray*}
	\int_{\S_{t_i} } \left(\frac{H_1H_{n-1}}{H_n}
	-1\right) \left(1+\cos\theta  \langle\nu,e\rangle\right)dA  \to 0, \hbox{ as }t_i\to +\infty.
	\end{eqnarray*}
Moreover one can show that for any sequence $t_i\to \infty$, there exists a convergent subsequence, whose
limit satisfying $$\left(\frac{H_1H_{n-1}}{H_n}
-1\right) \left(1+\cos\theta  \langle\nu,e\rangle\right)=0.$$
It is easy to see that the limit is a 
 spherical cap. 
	Next we show that any limit of a convergent subsequence  is uniquely determined, which implies the flow smoothly converges to a unique spherical cap. We shall use the argument in \cite{SWX}.
	
	Note that we have proved that $x(\cdot,t)$ subconverges smoothly to a capillary boundary spherical cap $ C_{\rho_\infty,\theta}(e_\infty)$.  Since $\V_{n,\theta}$ is preserved along  flow \eqref{flow with capillary}, the radius $\rho_\infty$ is independent of the choice of the subsequence of $t$.
	We now show in the following that $e_\infty=e$.  {Denote $\rho(\cdot,t)$ be the radius of the unique spherical cap $C_{\rho(\cdot, t), \theta}(e)$ around $e$ with contact angle $\theta$ passing through the point $x(\cdot,t)$.} Due to the
	spherical barrier estimate, i.e. Proposition \ref{barrier}, we know
	\begin{eqnarray*}
	\rho_{\max} (t):=\max \rho(\cdot,t)=\rho(\xi_t,t),
	\end{eqnarray*} is non-increasing with respect to $t$,  for some point $\xi_t\in M$.
Hence the limit $\lim\limits_{t\to +\infty} \rho_{\text{max}}(t)$ exists. Next we claim that
\begin{eqnarray}\label{max limit}
\lim_{t\to +\infty} \rho_{\max}(t)=\rho_\infty.
\end{eqnarray}
We prove this claim by contradiction. Suppose \eqref{max limit} is not true, then there exists $\varepsilon>0$ such that
\begin{eqnarray}\label{lim of radius}
\rho_{\max}(t)> \rho_\infty +{\varepsilon}, \hbox{ for } t \hbox{ large enough.}
\end{eqnarray} 
{By definition, $\rho(\cdot,t)$ satisfies
\begin{eqnarray}\label{spherical cap}
\rho^2\sin^2\theta=|x|^2-2\rho\cos\theta \langle x,e\rangle.
\end{eqnarray}}
Hence
\begin{eqnarray*}
 \left(\rho\sin^2\theta+\cos\theta \langle x,e\rangle \right) \p_t \rho
 =  \< \p_t x, x-\rho\cos\theta e\>.
\end{eqnarray*}
We evaluate at $(\xi_t,t)$. Since $\S_t$ is tangential to $C_{\rho,\theta }(e)$ at $x(\xi_t, t)$,
we have
$$\nu_{\S_t}(\xi_t, t)=\nu_{\p C_{ r,\theta}(e)}(\xi_t, t)= {\frac{x-\rho\cos\theta e}{\rho}}.$$
Thus we deduce
\begin{eqnarray}\label{r-max}
 {\left(\rho_{\max}\sin^2\theta+\cos\theta \langle x,e\rangle\right)\p_t\rho|_{(\xi_t,t)}}
= \rho_{\max}\left(\frac{1+\cos\theta  \langle\nu, e\rangle }{F}-\langle x,\nu\rangle\right). 
\end{eqnarray}
We note that there exists some $\delta>0$ such that \begin{eqnarray}\label{claim1}
{\rho_{\max}\sin^2\theta+\cos\theta \langle x,e\rangle}\ge \delta>0.
\end{eqnarray}In fact, this follows directly from \eqref{spherical cap}, due to
\begin{eqnarray}\label{spherical cap1}
{\rho \sin^2\theta+\cos\theta \langle x,e\rangle}=\frac{1}{2\rho}(|x|^2+\rho^2\sin^2\theta)\geq \frac{\rho}{2}\sin^2\theta>0.
\end{eqnarray}
Since the spherical caps $C_{\rho_{\max},\theta }(e)$ are the static solutions to \eqref{flow with capillary} and $x(\cdot,t)$ is tangential to $C_{\rho_{\max},\theta }(e)$ at $x(\xi_t,t)$, we see from \eqref{static model radius}
\begin{eqnarray}\label{xeq44}
\frac{ 1+\cos\theta\langle\nu,e\rangle}{\langle x,\nu\rangle }\Big|_{x(\xi_t,t)}= \frac{ 1+\cos\theta\langle\nu,e\rangle}{\langle x,\nu\rangle }\Big|_{C_{\rho_{\max},\theta }(e)}=\frac{1}{\rho_{\max}(t)}.
\end{eqnarray}
Since $x(\cdot,t)$ {subconverges} to $C_{\rho_\infty ,\theta }(e_\infty)$ and $\rho_\infty$ is uniquely determined, we have
\begin{eqnarray*}
F=\frac{n\sigma_n}{\sigma_{n-1}}\to \frac{1}{\rho_\infty} \quad \text{uniformly},
\end{eqnarray*}as $t\to+\infty$. 
Thus there exists $T_0>0$ such that 
\begin{eqnarray*}
\frac{1}{F} - {\rho_\infty} <\frac{\ep}{2},
\end{eqnarray*}
and hence 
\begin{eqnarray*}
\frac{1}{F} - {\rho_{\max}(t)} <-\frac{\ep}{2},
\end{eqnarray*}for all $t>T_0$.
Taking into account of \eqref{xeq44}, we see
\begin{eqnarray}\label{xeq45}
\left(\frac{1}{F} - {\frac{\langle x,\nu\rangle }{ 1+\cos\theta\langle\nu,e\rangle}}\right)\Big|_{x(\xi_t,t)} <-\frac{\ep}{2},
\end{eqnarray}for all $t>T_0$. {By adopting  Hamilton's trick, we conclude from \eqref{r-max}, \eqref{claim1} and \eqref{xeq45} that there exists some $C>0$ such that for almost every $t$, $$\frac{d}{dt}\rho_{\max}\le -C\ep.$$}
This is a contradiction  to the fact  that $\lim\limits_{t\to +\infty}\frac{d}{dt} \rho_{\max}=0$, and hence  claim \eqref{max limit} is true. Similarly, we can obtain that \begin{eqnarray}\label{min limit}
\lim_{t\to +\infty} \rho_{\min}(t)=\rho_\infty.
\end{eqnarray}
Hence  $\lim\limits_{t\to\infty} \rho(\cdot,t)=\rho_\infty$.  This implies that  any limit of a convergent subsequence is the spherical cap around $e$ with radius $\rho_{\infty}$. We complete the proof of Proposition \ref{converge limit}.

	\end{proof}
 In view of Proposition \ref{global existence} and Proposition \ref{converge limit}, Theorem \ref{longtime exist and conv} are proved.

\section{Alexandrov-Fenchel inequalities}\label{sect 5}

 In this section, we apply the convergence result of flow \eqref{flow with capillary} to prove Theorem \ref{thm 1}. 

\begin{proof}[\textbf{Proof of Theorem \ref{thm 1}}] 
Remember 	\begin{eqnarray}\label{AFI}
	\V_{k,\theta}(\widehat{C_{ r,\theta} })= r^{n+1-k} \bt,
	\end{eqnarray}
where $\bt$ was defined by \eqref{bt}.

	Assume that $\S$ is strictly convex. We have proved in Section 4 that   flow \eqref{flow with capillary} converges  a spherical cap, which we denote by  $C_{r_\infty,\theta}(e)$.
	 By the monotonicity of  $\V_{n,\theta}$ and  $\V_{k,\theta}$, Proposition \ref{monotone along flow} we have 
$$ \V_{n,\theta} (\widehat\Sigma) =  \V_{n,\theta} (\widehat{C_{r_\infty,\theta }(e)} ), \qquad
\V_{k,\theta} (\widehat\Sigma) \le   \V_{k,\theta} (\widehat{C_{r_\infty,\theta }(e)} ),$$
moreover, equality holds iff $\S$ is a spherical cap. It is clear that \eqref{AFI} is the same as
\eqref{af ineq}.

	When $\S$ is convex but not strictly convex, the inequality follows by approximation. The equality characterization can be proved  similar to \cite[Section 4]{SWX}, by using an argument of \cite{GL09}. We omit the details here.
\end{proof}

 For Corollary \ref{cor 1}, one just notes that when  $n=2$,
\begin{eqnarray*}
\bt=\frac{1}{3} \left(2-3\cos\theta +\cos^3\theta\right)\pi.
\end{eqnarray*} 

\
\

\section{
	Statements and Declarations.}

\

{\bf Conflict of interest.} On behalf of all authors, the corresponding author states that there is no conflict of
interest.

Data sharing not applicable to this article as no datasets were generated or analysed during the current study.

The authors have no relevant financial or non-financial interests to disclose.

\

\

\

\noindent\textbf{Acknowledgment:} LW is supported by NSFC (Grant No.  12201003, 12171260). CX is supported by NSFC (Grant No. 11871406, 12271449). We would like to thank the referee for   careful  reading and valuable suggestions to improve the context of the paper.

\end{document}